\documentstyle[12pt]{article}

\newcommand{\numero}[1]{
\addtocounter{section}{1}
\begin{center}{\bf \thesection .\
#1\vspace{-.1in}}\end{center}
\setcounter{subsection}{0}
\setcounter{lemma}{0}\indent}

\newcommand{\subnumero}[1]{
\pagebreak[1]\begin{center}{\em #1}\nopagebreak\end{center}
}

\newcommand{\eop}{\hfill $/$\hspace*{-.1cm}$/$\hspace*{-.1cm}$/$\vspace{.1in}}

\newtheorem{lemma}{Lemma}[section]
\newtheorem{theorem}[lemma]{Theorem}

\newtheorem{conjecture}{Conjecture}
\newtheorem{proposition}[lemma]{Proposition}

\newcommand{\cc}{{\bf C}}
\newcommand{\rr}{{\bf R}}
\newcommand{\qq}{{\bf Q}}

\newcommand{\pp}{{\bf P}}

\newcommand{\Cc}{{\cal C}}
\newcommand{\Ee}{{\cal E}}
\newcommand{\Ff}{{\cal F}}
\newcommand{\Gg}{{\cal G}}

\newcommand{\Oo}{{\cal O}}

\newcommand{\Aa}{{\cal A}}

\newcommand{\Mm}{{\cal M}}

\newcommand{\Ll}{{\cal L}}

\newcommand{\Gm}{{\bf G}_m}

\newcommand{\germ}{{\bf m}}

\newcommand{\delbar}{\overline{\partial}}

\begin{document}

\section*{Mixed twistor structures}

Carlos Simpson
\newline
CNRS, UMR 5580, Universit\'e Paul Sabatier\newline
31062 Toulouse CEDEX, France

\bigskip

The purpose of this paper is to introduce the notion of {\em mixed twistor
structure} as a generalization of the notion of mixed Hodge structure.  Recall
that a mixed Hodge structure is a vector space $V$ with three filtrations $F$,
$F'$ and $W$ (the first two decreasing, the last increasing) such that the two
filtrations $F$ and $F'$ induce $i$-opposed filtrations on $Gr^W_i(V)$.
(Generally $(V,W)$ is required to have a real structure and $F'$ is the
complex conjugate of $F$ but that is not very relevant for us here.)
Given a MHS $(V,W,F,F')$ we can form the {\em Rees bundle} $E:=\xi (V,F,F')$ over
$\pp ^1$ \cite{NAHT} \cite{SantaCruz}. In brief this is obtained from the
trivial
bundle $V\times \pp ^1$ by using $F$ to make an elementary transformation over
$0$ and $F'$ to make an elementary transformation over $\infty$.
The bundle $E$ is graded by strict subbundles which we denote $W_{i}E$ and
the condition of opposedness of the filtrations is equivalent to the condition
that $Gr^W_i(E)$ be semistable of slope $i$ on $\pp ^1$ (in other words a direct
sum of copies of $\Oo _{\pp ^1}(i)$).  The notion of mixed twistor structure is
simply obtained by abstracting this situation: an MTS is a pair consisting of a
bundle $E$ over $\pp ^1$ and a filtration by strict subbundles $W_iE$ such that
the $Gr^W_i(E)$ are semistable of slope $i$.  In the construction starting with
a mixed Hodge structure the resulting $(E,W)$ has an action of $\Gm$ covering
the standard action on $\pp ^1$ and in fact the mixed Hodge structures are
simply the $\Gm$-equivariant mixed twistor structures. Thus, in some sense, the
passage from ``Hodge'' to ``Twistor'' is simply forgetting to have an action of
$\Gm$.  This principle occured already, in a primitive way, in the passage from
systems of Hodge bundles (cf \cite{YMTU}) to Higgs bundles in \cite{CVHS}.

We will give some generalizations of basic classical results for
mixed Hodge structures, to the mixed twistor setting.  The process of
making these
generalizations is relatively direct although some work must be done to develop
the appropriate notion of {\em variation of mixed twistor structure}. The
overall
idea is that we have the following

\noindent
{\bf Meta-theorem} {\em
If the words ``mixed Hodge structure'' (resp. ``variation of mixed Hodge
structure'') are replaced by the words ``mixed twistor structure'' (resp.
``variation of mixed twistor structure'') in the hypotheses and conclusions of
any theorem in Hodge theory, then one obtains a true statement. The proof of the
new statement will be analogous to the proof of the old statement.}

We don't prove this meta-theorem but support it with several examples using the
basic theorems of mixed Hodge theory.

The utility
of the notion of mixed twistor structure, and of the above meta-theorem, is to
make possible a theory of weights for various things surrounding arbitrary
representations of the fundamental group of a smooth projective variety,
where up
until now the theory of weights has only been available for variations of Hodge
structure. This phrase needs the further explanation that a harmonic bundle
(cf \cite{Hitchin1} \cite{Corlette} \cite{HBLS}) yields a variation of pure
twistor structure (i.e. of mixed twistor structure with only one nonzero
weight-graded quotient) and in fact, up to choosing the weight a  variation of
pure twistor structure is essentially the same thing as a harmonic bundle. In
particular, any irreducible representation of $\pi _1(X,x)$ for a compact
K\"ahler manifold $X$, underlies a variation of pure twistor structure unique up
to shift of weight (i.e. tensorization with constant pure rank one twistor
structures). This is explained in Lemma \ref{harmonic} below.

Apart from the above meta-theorem as a way of obtaining new statements, a
possible
area where we obtain a new type of object is the following: over the moduli
space
$M_B$ of representations of $\pi _1(X)$ for a compact K\"ahler variety $X$, we
may have several families of mixed twistor structures, for example cohomology of
open or singular subvarieties of $X$ with coefficients in the VMTS corresponding
to $\rho \in M_B$, or some types of rational homotopy invariants relative to the
representation $\rho$ (such as the relative Malcev completion or its analogues
for higher homotopy).  These natural families of MTS give classifying maps
$M_B\rightarrow {\cal MTS}$ to the moduli stack of mixed twistor structures
(the maps will not usually be algebraic but should be real analytic, for
example).  It might be interesting to study these classifying maps.

The terminology ``twistor'' comes from the particular case of weight one twistor
structures with a certain kind of real structure ({\em antipodal}).  These are
equivalent to the twistor bundles over $\pp ^1$ of quaternionic vector spaces,
see \cite{HKLR} for a nice exposition.
Deligne in \cite{DeligneLetter} originally explained to me how to associate a
quaternionic vector space to a weight one real Hodge structure and how to
interpret the twistor space for Hitchin's hyperk\"ahler structure in terms
of moduli of $\lambda$-connections, a notion very closely related to the Hodge
filtration on the de Rham complex. This eventually led to a study of the formal
neighborhood of the twistor lines in the twistor space, to an interpretation
in terms of mixed Hodge structures (for the case of the moduli space of
representations), and consequently to the present definition. In \S 7 below we
discuss the relation between mixed twistor structures and the formal power
series for hyperk\"ahler (or more precisely any hypercomplex) structures at a
point. The case of the moduli space of
representations of $\pi _1$ of a smooth projective variety is taken up
briefly in
\S 8.

Here is an outline of the paper.

{\em \S 1---Mixed twistor structures:}
definition of mixed twistor structures, and the relation with mixed Hodge
structures. The theorem that the category of MTS is abelian. As an aside at the
end we point out that this theorem works with $\pp ^1$ replaced by any
projective variety as base (this remark is not used later, though).

{\em \S 2---Real structures:}
two different kinds of real structures for a mixed twistor structure: circular
and antipodal. The weight one pure antipodal real twistor structures are the
same as the quaternionic vector spaces.

{\em \S 3---Variations of mixed twistor structure:}
definition of  $\Cc ^{\infty}$ families of
mixed twistor structure, variations of mixed twistor structure, polarizations,
real structures and a holomorphic interpretation as a pair of
$\lambda$-connections on the two standard affine lines in $\pp ^1$.

{\em \S 4---Cohomology
of smooth compact K\"{a}hler manifolds with VMTS coefficients:}
we treat the cohomology of smooth compact K\"ahler manifolds, first with
variations of pure twistor structure as coefficients, then with mixed
coefficients. We give a holomorphic interpretation using $\lambda$-connections.

{\em \S 5---Cohomology of open and singular varieties:}
generalization of the basic theorems about cohomology of open and singular
varieties to the mixed twistor case. This involves a notion of {\em mixed
twistor
complex}. A somewhat new point here involves the notion of {\em patching}
complexes of sheaves on open sets of $\pp ^1$; this replaces certain arguments
involving the real structure in the theory of mixed Hodge complexes.

{\em \S 6---Nilpotent orbits and the limiting mixed twistor structure:}
we give a conjectural version of the nilpotent orbit theorem relating
the degeneration of a harmonic bundle on a punctured disc to a limiting mixed
twistor structure.

{\em \S 7---Jet bundles of hypercomplex manifolds:}
description of the formal germ of a hypercomplex structure in the
neighborhood of
a point by a mixed twistor structure on the jet space at the point.  This comes
from looking at the jet bundle of the twistor space along a twistor line.
Following the philosophy described above of looking at classifying maps, we
define
the {\em Gauss map} from a hypercomplex manifold into the moduli stack for mixed
twistor structures on the jet spaces.

{\em \S 8---The moduli space of representations:}
this is a discussion similar to that of the previous section, for moduli spaces
of representations of fundamental groups of smooth projective varieties.
Using the family of moduli spaces $M_{Del}$ of \cite{SantaCruz} as a replacement
for the twistor space, we are able to treat singular points: we get a mixed
twistor structure on the generalization of the jet space at a singular point,
and this determines the formal neighborhood of a twistor line.

{\em Notations:}
we work in the analytic category of complex analytic spaces (usually smooth
complex varieties) and when talking about projective varieties we mean the
associated analytic varieties with the usual topology.

\numero{Mixed twistor structures}

Fix the projective line $\pp ^1$ with points $0,1,\infty$, with the standard
line bundle $\Oo _{\pp ^1}(1)$.

A {\em twistor structure} is a vector bundle $E$ on $\pp ^1$.  The {\em
underlying vector space} is $E_1$, the fiber of the bundle $E$ over $1\in \pp
^1$.  We say that a twistor structure $E$ is {\em pure of weight $w$} if $E$ is
a semistable vector bundle of slope $w$, which is equivalent to saying that $E$
is a direct sum of copies of $\Oo _{\pp ^1}(w)$.

A {\em mixed twistor structure} is a twistor structure $E$ filtered by an
increasing sequence of strict subbundles $W_iE$ such that for all
$i$, $Gr^W_i(E)=
W_iE / W_{i-1}E$ is pure of weight $i$. The filtration $W_{\cdot}E$ is called
the {\em weight filtration}.  A mixed twistor structure is said to be {\em pure}
if the associated graded is nontrivial in only one degree.

\subnumero{Relation with mixed Hodge structures}

As described in \cite{NAHT} and (\cite{SantaCruz} \S 5), if $V$ is a
complex vector space with two decreasing filtrations $F$ and $F'$, then we
obtain
a vector bundle $\xi (V; F, F')$ on $\pp ^1$. This is the {\em twistor structure
associated to the pair of filtrations $F$ and $F'$}.  The two filtrations define
a pure Hodge structure of weight $w$ if and only if $\xi (V; F, F')$
is a pure twistor structure of weight $w$.

Suppose $V$ is a vector space with three filtrations $W_{\cdot}$ (increasing)
and $F$ and $F'$ (decreasing). Then putting
$$
W_i \xi (V; F,F'):= \xi (W_iV; F, F')
$$
we obtain a filtration by strict subbundles. The associated-graded
$Gr ^W_i(\xi (V; F,F')$ is pure of weight $w$ if and only if $F$ and $F'$ induce
a pure Hodge structure of weight $w$ on $Gr ^W_i(V)$.

The condition in the previous paragraph is what we call a {\em complex mixed
Hodge structure} (a definition which should have been made long ago).  Note
that a
real mixed Hodge structure is simply a real vector space $V_{\rr}$ with
filtration
$W_{\cdot}$ and complex filtration $F$ of  $V_{\cc}$ such that $W, F,
\overline{F}$ form a complex mixed Hodge structure on $V_{\cc}$.
Similarly if $V,W,F,F'$ is a complex MHS then $V\oplus \overline{V}$ has a
structure of real MHS (using the real structure which interchanges $V$ and
$\overline{V}$).

We have the following characterization.

\begin{lemma}
If $V_{\rr}$ is a real vector space with increasing real filtration $W$, and
with complex filtration (decreasing) $F$ of $V_{\cc}$, then these data define a
real mixed Hodge structure if and only if $\xi (V_{\cc}; F, \overline{F})$ with
the induced filtration $W_{\cdot}\xi (V_{\cc}; F, \overline{F})$ is a mixed
twistor structure.
\end{lemma}
\eop

The group $\Gm$ acts on $\pp ^1$ by translation. Denote the morphism of this
action by
$$
\mu : \Gm \times \pp ^1 \rightarrow \pp ^1.
$$
A {\em $\Gm$-equivariant mixed twistor
structure} is a mixed twistor structure $(E, W)$ together with an isomorphism
of filtered coherent sheaves on $\Gm \times \pp ^1$
$$
\rho : \mu ^{\ast}(E,W)  \cong p_2^{\ast}(E,W)
$$
such that the two resulting morphisms on $\Gm \times \Gm \times \pp ^1$
$$
(\mu \circ \mu )^{\ast}(E,W) \cong p_3^{\ast}(E, W)
$$
are equal.  If we view $E$ as a filtered vector bundle over $\pp ^1$ this just
means that the action of $\Gm$ is lifted to an action on the total space of the
bundle preserving the filtration.

A {\em morphism} of equivariant mixed twistor structures is a morphism of
filtered bundles compatible with the action.

\begin{proposition}
The category of $\Gm$-equivariant mixed twistor structures is naturally
equivalent to the category of complex mixed Hodge structures (see below).
\end{proposition}
{\em Proof:}
As was seen in \cite{NAHT} and \cite{SantaCruz} a bundle over $\pp ^1$
equivariant for the action of $\Gm$ is the same thing as a vector space $V=E_1$
together with two decreasing filtrations $F$ and $F'$: the  bundle $E$ is
recovered as  $\xi (V, F, F')$.

Restriction to the fiber $V= E_1$ over $1\in
\pp ^1$ induces an inclusion-preserving bijection between the set of strict
subsheaves of $V$ preserved by $\Gm$, and the set of subspaces of $V$. To see
this, the restriction to the fiber may be viewed as the composition of
restriction to $\Gm\subset \pp ^1$ then restriction to the fiber. The second
arrow is obviously a bijection between subsheaves of $E|_{\Gm}$ preserved by
$\Gm$ and subspaces of $E_1$ (using the action of $\Gm$). Note by the way
that any equivariant subsheaf of $E|_{\Gm}$ is automatically strict.  Recall
that a strict subsheaf of a vector bundle over a curve is determined by its
restriction to any Zariski open set; this implies that restriction is a
bijection
between all strict subsheaves of $E$ and all strict subsheaves of
$E|_{\Gm}$ (and
this bijection preserves the subset of those which are equivariant).

In view of this bijection, we obtain a bijection between the set of
increasing
filtrations of $E$ by strict subsheaves preserved by $\Gm$, and the set of
increasing filtrations of $V$ by subspaces.  In view of this and  the first
paragraph of the proof, an equivariant  vector bundle $E$ over $\pp ^1$
with increasing filtration by strict subsheaves $W_nE$ is the same thing as a
vector space $V=E_1$ together with an increasing filtration $W_nV$ and
two decreasing filtrations $F$ and $F'$.  The construction in the reverse
sense is the Rees bundle construction $E = \xi (V, F, F')$. Finally note that
$$
Gr ^{WE}_n(E) = \xi ( Gr ^{WV}_n(V), F_{Gr}, F'_{Gr})
$$
where $F_{Gr}$ and $F'_{Gr}$ are the filtrations induced on the
associated-graded (the Rees-bundle construction is compatible
with taking strict subquotients, as can be verified for each filtration $F$
and $F'$ separately using splittings).  The two filtrations
$F_{Gr}$ and $F'_{Gr}$ are $n$-opposed if
and only if the Rees bundle is semistable of slope $n$.  Thus the data $(V,W,F,
F')$ defines a mixed Hodge structure if and only if $(E,W)$ is a mixed twistor
structure.

We leave to the reader to verify full faithfulness of this correspondence
for morphisms.
\eop

We don't actually use the above proposition anywhere but it is crucial for
understanding the analogy between all of the results which we give below, and
the corresponding results for mixed Hodge structures. A somewhat pertinent
remark is that our proofs for mixed twistor structures are all $\Gm$-equivariant
in case the input-data is $\Gm$-equivariant. Thus our proofs give proofs of the
corresponding statements for mixed Hodge structures.  One hesitates to say that
these constitute ``new'' proofs since they are really just recopying the old
proofs for mixed Hodge structures into this new language.

\subnumero{The abelian category of mixed twistor structures}

The first basic theorem in the theory of mixed Hodge structures is the fact that
the  category of MHS is abelian. This generalizes to mixed twistor structures
and the proof is essentially the same.

\begin{lemma}
\label{abelian}
{\rm (cf \cite{Hodge2}, 2.3.5(i))}
The category of mixed twistor structures is abelian.
\end{lemma}
{\em Proof:}
Suppose $f: (E, W) \rightarrow (E', W')$ is a morphism of mixed twistor
structures.

{\bf Step 1}\,\,
Define the {\em cokernel} of $f$ to be $(E'', W'')$ where
$E''$ is the cokernel of $f$ considered as a map of coherent sheaves, and $W''$
is the filtration of $E''$ induced by $W'$.   For the moment, this is just a
coherent sheaf filtered by subsheaves. Let $A \subset E$ denote the kernel
subsheaf of $f$.

We work by induction on the size of the interval where the two weight
filtrations $W$ and $W'$ are supported.  If this interval has size $1$ then
the cokernel is just the cokernel of a map of semistable sheaves of the same
slope, so the cokernel is also semistable.  Now consider an interval size of at
least two.  Let $Gr ^W_i(E)$ and $Gr ^{W'}_i(E')$ be the highest nonzero pair of
terms (in other words $i$ is the smallest integer such that $W_i=E$ and $W'_i =
E'$). Then  we have an exact sequence $$
Gr ^W_i(E) \rightarrow Gr ^{W'}_i(E')
\rightarrow Gr ^{W''}_i(E'') \rightarrow 0.
$$
In particular, $Gr ^{W''}_i(E'')$ is semistable of slope $i$.  Let $V_{i-1}
\subset E$ be the inverse image of $W'_{i-1}$. We have an exact sequence
$$
0 \rightarrow W_{i-1} \rightarrow V_{i-1} \rightarrow B\rightarrow 0
$$
where $B$ is the kernel of the map
$$
Gr ^W_i(E) \rightarrow Gr ^{W'}_i(E') .
$$
Since this is a map of semistable sheaves of slope $i$, the kernel $B$ is
semistable of slope $i$.  Thus $V_{i-1}$ is a mixed twistor structure (with the
weight filtration induced by that of $E$).

By the induction hypothesis, the cokernel of the morphism
$W_{i-1}\rightarrow W'_{i-1}$ of mixed twistor structures concentrated in a
smaller interval (the $i$-th graded pieces are zero this time), is again a mixed
twistor structure concentrated in the same interval.  In particular, $W'_{i-1}
/f(W_{i-1})$ cannot support a morphism from a semistable sheaf of weight
$\geq i$.
Thus the morphism
$$
B = V_{i-1} /W_{i-1} \rightarrow W'_{i-1} /f(W_{i-1})
$$
is zero.   Thus
$$
f(V_{i-1}) = f(W_{i-1}).
$$
Note that
$$
f(V_i) = f(E) \cap W'_{i-1} = \ker (W'_{i-1} \rightarrow W''_{i-1}).
$$
Therefore we get an exact sequence
$$
W_{i-1} \rightarrow W'_{i-1} \rightarrow W'' _{i-1} \rightarrow 0.
$$
Furthermore the weight filtration on $W'' _{i-1}$ is that induced by the weight
filtration of $W'_{i-1}$.  Hence $W'' _{i-1}$ is the cokernel of the map
$W_{i-1} \rightarrow W'_{i-1}$, and by our induction hypothesis this cokernel is
a mixed twistor structure.  Finally, $E''$ is an extension of $Gr ^{W''}_i(E'')$
(which we have seen to be pure of weight $i$ at the start) by $W'' _{i-1}$, thus
$E''$ with its weight filtration $W''$ is a mixed twistor structure. Along the
way, we have proved that $E''$ is a bundle, for $Gr ^{W''}_i(E'')$ is a bundle
(since it is a cokernel of a map of semistable bundles of the same slope)
and $W''
_{i-1}$ is a bundle by induction.

{\bf Step 2}\, \,
Let $A \subset E$ be the kernel of $f$ considered as a morphism
of coherent sheaves, and put $W_i A := A \cap W_i$. Note that $A$ is a strict
subbundle of $E$ and $W_iA$ are strict subbundles of $A$.  We would like
to show that $(A, W_{\cdot }A)$ is a mixed twistor structure, i.e. to
show that $Gr ^{WA}_i(A)$ is semistable of slope $i$.

(Actually the following proof is just the dual of the proof of step 1.)

Let $i$ denote the smallest number where either $W_i$ or $W'_i$ is nonzero.
Thus $W_i$ and $W'_i$ are both semistable bundles of slope $i$.

Let $V_i := f^{-1}(W'_i)\subset E$.

We have that $V_i/W_i$ is the kernel of the morphism of mixed twistor
structures $E/W_i \rightarrow E'/W'_i$. By induction (on the size of the
interval where the weight graded quotients of $E$ and $E'$ live, just as in Step
1) we can suppose that $V_i/W_i$ with its filtration induced by $W_j/W_i$ is a
mixed twistor structure, substructure of $E/W_i$. In particular the weights of
$V_i/W_i$ are $> i$, which implies that
any morphism of $V_i / W_i $ to a semistable bundle of slope $\leq i$ is
zero. Thus the morphism $V_i/W_i \rightarrow W'_i / f(W_i)$ is zero.
This implies that
$$
V_i = W_i + A,
$$
and from this formula the sequence
$$
0\rightarrow A / W_{i}A \rightarrow  E / W_{i}
\rightarrow E' / W' _{i} \rightarrow 0
$$
is exact. Note also that $A/W_iA = V_i /W_i$.

For $j \geq i$ we have $W_jA / W_iA = A / W_iA \cap W_j /W_i$.
In particular the filtration induced on $A/W_iA$ by the filtration $W_iA$ is the
same as the kernel filtration for the above exact sequence.  By induction on the
size of the range where the weight filtration lives, we have that $A/W_iA$ is a
mixed twistor structure.

On the other hand, $W_iA $ is the kernel of the map $W_i \rightarrow W'_i$
of semistable bundles of slope $i$, so $W_iA $ is semistable of slope $i$.
Note of course that $W_{i-1}A=0$. Thus for any $j$ we have that $Gr ^{WA}_j(A)$
is semistable of slope $j$ and $A$ is a mixed twistor structure.

{\bf Step 3}\,\,
We have to show that the image is equal to the coimage. In view of what we have
seen with cokernels and kernels, this comes down to showing that if $f$ is an
isomorphism on coherent sheaves $E \cong E'$ then $f(W_i)= W'_i$.

Again, suppose that this is true for pairs of MTS with weight filtrations
concentrated in any smaller range. Let $W_i$ and $W'_i$ be the lowest nonzero
pair of levels in the weight filtration.  Above we established in this case the
formula
$$
f^{-1}(W'_i ) = W_i + {\rm ker}(f)
$$
but by hypothesis ${\rm ker}(f)=0$ here so $f^{-1}(W'_i )= W_i$.
Since $f$ is an isomorphism, $f(W_i)=W'_i$.  Now we can look at the isomorphism
$E/W_i \cong E' /W'_i$ and by our inductive hypothesis the weight filtrations
there are the same.  Thus $f$ induces for all $j$ an isomorphism $Gr
^W_j(E)\cong Gr ^{W'}_j(E')$, and this implies that $f(W_j)=W'_j$.

We have now completed all of the verifications necessary to show that our
category is abelian (cf \cite{Hodge2}).
\eop

{\em Remark:}  The reader may generalize the other parts (ii)--(v) of
Theorem 2.3.5 of \cite{Hodge2}. Note that the correct generalization of part (v)
is
that the functor ``fiber over $0\in \pp ^1$'' is exact, and this may be further
generalized to the statement that the functor ``fiber over $t\in \pp ^1$''
is exact for any point $t$ (in the Hodge case all points are equivalent to $0$,
$1$ or $\infty$ by the $\Gm$ action; the point $\infty$ corresponds to the
complex conjugate filtration, so nothing is new in the Hodge case but in the
twistor case all points are different).

\subnumero{Aside: some generalizations of the above theorem}
The above theorem that we get an abelian category depends only on some fairly
standard results about semistability. These results hold in much greater
generality than just semistable bundles on $\pp ^1$.  We could make the
following definition. If $X$ is a projective scheme then a {\em mixed sheaf of
pure dimension $d$ on $X$} is a sheaf $E$ of pure dimension $d$ with a
filtration
$W_P$ by subsheaves indexed by polynomials $P\in \qq [x]$ (increasing for the
ordering of polynomials by values at large $x$) such that there are only a
finite number of jumps in the filtration and such that  $Gr _P^W(E)$ is
a semistable sheaf of pure dimension $d$ with normalized Hilbert polynomial $P$
\cite{Moduli}.  By the same proof as above, the category of mixed sheaves of
pure dimension $d$ is abelian.  When $d=dim (X)$ these become torsion-free
sheaves.

We can generalize further if anybody is interested!  Let $\Lambda$ be a sheaf
of rings of differential operators such as considered in \cite{Moduli}.
Then a {\em mixed $\Lambda$-module} is a $\Lambda$-module $E$ of some pure
dimension $d$ with a filtration by sub-$\Lambda$-modules $W_P$ indexed by
rational polynomials $P$, with a finite number of jumps, such that
$Gr ^W_P(E)$ is a semistable $\Lambda$-module with normalized Hilbert polynomial
$P$.

There are probably other variants, for example with parabolic structure, etc.

\subnumero{Moduli of mixed twistor structures}

There is an obvious notion of algebraic family of mixed twistor structures:
if $S$ is a scheme then a family of MTS over $S$ is a bundle $E$ on $S\times \pp
^1$ provided with a filtration by strict subbundles $W_iE$ such that in the
fiber
over each point $s\in S$ this gives a MTS.

If we associate to each scheme $S$ the category of families of MTS over
$S$, we obtain a stack ${\cal MTS}^{\rm any}$. The superscript denotes the
fact that we include any morphisms of MTS, in particular it is not a stack of
groupoids. The associated stack of groupoids which we denote by
${\cal MTS}$ is an algebraic stack in the sense of Artin. It is locally of
finite
type as we shall see below: the pieces corresponding to fixing the dimensions
of the weight-graded pieces, are of finite type.

 There is a
natural action of $\Gm$ and the fixed points (in the stack-theoretic sense) are
exactly the $\Gm$-equivariant mixed twistor structures, i.e. complex mixed Hodge
structures.

We will attempt briefly to give some idea of what ${\cal MTS}$ looks like.
A {\em framed $MTS$}  is a mixed twistor structure $(E,W)$ provided with
isomorphisms
$$
\beta _n: Gr ^W_n(E) \cong \Oo _{\pp ^1}(n)^{b_n}.
$$
This is equivalent to the data of frames for the vector spaces $Gr ^W_n (E_1)$.
A {\em framed family of MTS} over a base scheme $S$ is just a MTS over $S$
provided with isomorphisms $\beta _n$ on $\pp ^1 \times S$ as above.

Any automorphism of a framed MTS (or family of framed MTS) fixing the framing,
is the identity. To prove this, suppose that $f$ is such an automorphism.
Then $f - 1$ induces zero on the associated graded pieces, hence it is zero
(because kernel and cokernel commute with taking $Gr ^W$---but this is easy to
prove along the lines of the above theorem anyway).

We claim that the functor which to any scheme $S$ associates the set of
isomorphism classes of framed families of MTS over $S$, is representable by a
scheme $Fr {\cal MTS}$ locally of finite type (of finite type if we fix the
dimensions of the graded pieces i.e. the $b_n$ above). Let
$Fr {\cal MTS}(b_{\cdot})$ denote the part corresponding to MTS where $Gr
^W_n(E)$ has rank $b_n$. The group
$$
GL(b_{\cdot}):= \prod _{n} GL(b_n)
$$
acts on $Fr {\cal MTS}(b_{\cdot})$ with stack-theoretic quotient
${\cal MTS}(b_{\cdot})$ (the open and closed substack of ${\cal MTS}$
with given dimensions of graded pieces). The stack
${\cal MTS}$ is the disjoint union of the ${\cal MTS}(b_{\cdot})$
so it is locally of finite type.

We now prove the representability claim of the previous paragraph.  Write
$b_{\cdot } = (b_0, \ldots , b_k)$ with the rest being zero.  We may proceed by
induction on $k$, the case $k=0$ being obvious (the representing scheme is just
a point).  Thus we may assume that
$Y=Fr {\cal MTS}(0, b_1, \ldots , b_k)$ exists. There is a universal
family of framed MTS $(E', W')$ on $Y\times \pp ^1$. Suppose $S$ is a scheme and
$(E_S, W)$ is a family of framed MTS on $S$ with dimensions $b_0, \ldots , b_n$.
Then $E_S/W_1E_S$ is a family of framed MTS with dimensions $0, b_1, \ldots ,
b_k$ so we obtain a morphism $p:S\rightarrow Y$ and
$p^{\ast}(E',W')=(E_S/W_1,W)$
with framings. The data of $(E_S, W)$ is equivalent to the data
of $(E_S/W_1,W)$ together with an extension of the bundle
$E_S/W_1$ by $\Oo _{\pp ^1\times S}(0)^{b_0}$.  Thus we have a lifting
of $p$ to a morphism from $S$ into the relative $Ext^1$-bundle
(which we show to be a bundle two paragraphs below)
$$
Ext ^1_{\pp ^1\times Y/Y}(E', \Oo _{\pp ^1\times Y}^{b_0})\rightarrow Y.
$$
This lifting is unique, and the pullback of the universal extension gives
$(E_S, W)$.  Thus $Fr{\cal MTS}(b_0, \ldots , b_k)$ is the relative
$Ext^1$-bundle refered to above.

To sum up the previous paragraph, there are schemes
$Fr{\cal MTS}(b_j, \ldots , b_k)$ and universal families of framed MTS
$(E^{\rm univ}(b_j, \ldots , b_k), W)$ which are filtered bundles over
$$
Fr{\cal MTS}(b_j, \ldots , b_k)\times \pp ^1.
$$
There are natural morphisms
$$
Fr{\cal MTS}(b_j, \ldots , b_k)\rightarrow
Fr{\cal MTS}(b_{j+1}, \ldots , b_k)
$$
such that the upper scheme is naturally identified with the
 relative $Ext^1$-bundle classifying extensions in the $\pp ^1$-direction of
$E^{\rm univ}(b_{j+1}, \ldots , b_k)$ by $\Oo _{\pp ^1} (j)^{b_j}$.

We need to point out that the relative $Ext^1$ referred to above are indeed
bundles over the base $Y$. This is by semicontinuity theory using the fact that
for any point the corresponding $Ext^0$ is zero. Indeed the fiber over any point
$y\in Y$ of the bundle  $E'=E^{\rm univ}(b_{j+1}, \ldots , b_k)$ (over $\{ y\}
\times \pp ^1$) is an extension of bundles which are stable of slopes $\geq j+1$
so it decomposes into factors of degree $\geq j+1$. Thus there are no
homomorphisms to $\Oo _{\pp ^1} (j)^{b_j}$ so the $Ext^0=Hom$ is zero.
Since we work on $\pp ^1$ there are no other terms and semicontinuity implies
that the family of $Ext^1$ spaces all have the same dimension and fit together
into a vector bundle which we have referred to as
$Ext ^1_{\pp ^1\times Y/Y}(E', \Oo _{\pp ^1\times Y}^{b_0})$.

From this description, $Fr{\cal MTS}(b_j, \ldots , b_k)$ is obtained as
a sequence of vector bundles, eventually over a point. Thus the underlying
topological space is contractible. Hence the homotopy type of
${\cal MTS}(b_j, \ldots , b_k)^{\rm top}$ is the same as that of $BGL(b_j,
\ldots
, b_k)^{\rm top}$.

We give a formula for the dimension of the stack ${\cal MTS}(b_j, \ldots , b_k)$
(recall that this means roughly speaking the dimension of the space of orbits
minus the dimension of the stabilizer group). The dimension is
equal to
$$
dim (Fr{\cal MTS}(b_j, \ldots , b_k))- (b_j^2 + \ldots + b_k^2)
$$
the latter term  being the dimension of $GL(b_j, \ldots ,
b_k)$. For $dim(Fr{\cal MTS}(b_j, \ldots , b_k))$
we use the above expression as a series of fiber bundles.  Note that the
dimension of
$$
Ext ^1(\Oo _{\pp ^1}(i), \Oo _{\pp ^1}(j))= H^1(\pp ^1, \Oo _{\pp
^1}(j-i))
$$
is $i-j-1$ when $i>j$, by Riemann-Roch.  Thus we have
$$
dim (Fr{\cal MTS}(b_j, \ldots , b_k))=
dim (Fr{\cal MTS}(b_{j+1}, \ldots , b_k)) + \sum _{i=j+1}^{k}
(i-j-1)b_ib_j.
$$
Thus
$$
dim (Fr{\cal MTS}(b_j, \ldots , b_k))=\sum _{j\leq u < i \leq k}
(i-u-1)b_ib_j
$$
and the term $dim \, GL(b_j,\ldots , b_k)$ fits in nicely to give
$$
dim ({\cal MTS}(b_{\cdot}))=\sum _{u \leq i }
(i-u-1)b_ib_j .
$$
Note in particular that for most well spaced-out sequences $b_{\cdot}$ the
dimension is positive, and in particular the family of orbits is nontrivial.
This means that there are ``moduli'' of mixed twistor structures.

{\em Questions:}  Define a natural $GL(b_{\cdot})$-linearized line bundle
on $Fr{\cal MTS}(b_{\cdot})$. What are the stable or semistable points for the
action, and what does the GIT quotient space look like? Are the mixed twistor
structures which come up in nature (e.g.  the mixed Hodge structures)
semistable?

\numero{Real structures}
To complete the circle of definitions and comparisons, we define some notions of
{\em real twistor structure} and {\em real mixed twistor structure}.
The situation is more complicated than for mixed Hodge structures, because
there are many possible antiholomorphic involutions of $\pp ^1$.
We will isolate for our purposes two examples.

Let $\sigma _{\pp ^1}$ denote the antipodal involution of $\pp ^1$ (it is
antilinear, interchanges $0$ and $\infty$ and interchanges $1$ and $-1$).
An {\em antipodal real twistor structure} is a bundle $E$ on $\pp ^1$ with
antilinear involution $\sigma$ lying over $\sigma _{\pp ^1}$.
An {\em antipodal real mixed twistor structure} is an antipodal real twistor
structure with filtration by strict subbundles preserved by $\sigma$.

In considering polarizations in the next section, we will make use of the
following sheaf-theoretic point of view. If $E$ is a vector bundle on $\pp ^1$
considered as a locally free sheaf of $\Oo _{\pp ^1}$-modules then put
$$
\sigma ^{\ast}(E) (U):= E(\sigma _{\pp ^1}U).
$$
The structure of $\Oo _{\pp ^1}(U)$-module is defined by
setting, for $e\in E(\sigma _{\pp ^1}U)$ and $a\in \Oo_{\pp ^1}(U)$,
$a\cdot e := \overline{\sigma _{\ast}(a)}e$.  Note that
$\sigma ^{\ast}(E)$ is again a locally free sheaf of $\Oo _{\pp ^1}$-modules.
However the functor $\sigma ^{\ast}$ is antilinear: the induced map on
$\cc$-vector spaces
$$
Hom (E, E')\rightarrow Hom (\sigma ^{\ast}(E),\sigma ^{\ast}(E'))
$$
is antilinear.  Also there is an induced map
$$
H^i(E) \rightarrow H^i(\sigma ^{\ast}(E)),
$$
again antilinear.
Finally the functor $\sigma ^{\ast}$ is an involution (its square is naturally
equal to the identity).  An involution of $E$ lying over $\sigma _{\pp ^1}$
(the notion defined in the previous paragraph, which concerns the total space
of the bundle) can be expressed in sheaf-theoretic language as a morphism
$$
f: E \rightarrow \sigma ^{\ast}(E)
$$
such that $\sigma ^{\ast}(f)\circ f = 1_E$.

Let $\tau _{\pp ^1}$ denote the antiholomorphic involution which preserves the
unit circle.  A   {\em circular real twistor structure} is a bundle $E$ on $\pp
^1$ with antilinear involution $\tau$ lying over $\tau _{\pp ^1}$.
A {\em circular real mixed twistor structure} is a circular real twistor
structure
with filtration by strict subbundles preserved by $\tau$. Note that a
sheaf-theoretic discussion similar to that which we have given for $\sigma$,
exists for $\tau$ also.

The difference between these two notions is that $\sigma _{\pp ^1}$ has no fixed
points, whereas $1$ is a fixed point of $\tau _{\pp ^1}$. Thus if $E$ has a
circular real structure then $E_1$ is a real vector space, whereas this is not
necessarily the case for an antipodal real structure.

We now investigate more closely what these structures mean in the pure case.
There are two different situations depending on the parity of the weight.

Suppose $E$ is an antipodal (resp. circular) real twistor structure which
is pure
of weight $0$. Then  $H^0(\pp ^1, E)$ is a complex vector space with antilinear
involution induced by $\sigma$ (resp. $\tau$), or in other words it is a real
vector space. Note that
$$
E = H^0(\pp ^1, E)\otimes _{\cc } \Oo _{\pp ^1}
$$
preserving the involutions.
Conversely any
real vector space gives rise to an antipodal (resp. circular) real twistor
structure which is pure of weight $0$ and these constructions are inverses.
In the circular case note that the isomorphism
$$
H^0(\pp ^1, E) \cong E_1
$$
is compatible with the involution $\tau$.

Suppose $E$ is an antipodal real twistor structure which is pure of weight
$1$. Then we obtain a structure of quaternionic vector space.
Conversely a quaternionic vector space gives an antipodal real twistor structure
pure of weight $1$, and these constructions are inverses. These are both well
described in \cite{HKLR}. This is the example which lends the name
``twistor''. We
start by describing this second direction which is just the classical ``twistor
space'' construction of Penrose (cf \cite{Penrose} \cite{Hitchin2}
\cite{HKLR}).

A quaternionic vector space is a real vector space $V$ with operations $I$, $J$
and $K$ satisfying the relations
$$
I^2=J^2=K^2 = -1, \;\;\; IJ=K \;\; \mbox{etc.}
$$
For any triple of reals $(x,y,z)$ with $x^2+ y^2 + z^2= 1$ we can define a
complex structure $xI + yJ + zK$ on $V$. Combined with the standard complex
structure on $S^2= \pp ^1$ this serves to define an almost complex structure on
the bundle $E=V\times \pp ^1$.  This almost complex structure is integrable, as
can be checked by a direct calculation, or else by noting that the bundle
$L:=\Oo
_{\pp ^1}(1) \oplus \Oo _{\pp ^1}(1)$ has an antipodal real involution
wherefrom the invariant sections give a trivialization of the real bundle
$L\cong \pp ^1 \times {\bf H}$ such that the points $0$, $1$ and
$i$ in $\pp ^1$ give a quaternionic triple of complex structures $I$, $J$
and $K$
on ${\bf H}$ whose resulting twistor almost complex structure recovers $L$
(this
is a calculation which needs to be done and for which we refer to  \cite{HKLR}).
As any quaternionic $V$ is a direct sum of copies of ${\bf H}$ we obtain the
integrability of the twistor bundle $E$ (which is a direct sum of copies of
$L$),
as well as the fact that the twistor bundle is pure of weight $1$ with antipodal
real structure.

We now describe how an antipodal real twistor structure $E$ pure of weight $1$
comes from a quaternionic vector space. The underlying
real vector space $A$ is the space of sections of $E$ which are preserved by the
involution $\sigma$.  Note that $\sigma$ induces an antilinear involution on
$H^0(\pp^1, E)$ and the fixed points form a real subspace
$$
A= H^0(\pp^1, E)^{\sigma }
$$
whose dimension is
half the complex dimension of  $H^0(\pp^1, E)$. In turn, the complex
dimension of
$H^0(\pp^1, E)$ is twice the complex rank of $E$ (by purity of weight $1$) so
$$
dim _{\rr}(A)= rk (E).
$$
This suggests that the evaluation morphism $A\rightarrow E_p$ for any point
$p\in \pp ^1$, should be an isomorphism. We prove this: if $e\in E_p$ then
$\sigma (e)\in E_{\sigma p}$ and by  purity of weight $1$ there is a unique
section $f: \pp ^1 \rightarrow E$ such that $f(p)=e$ and $f(\sigma p)= \sigma
(e)$.  Uniqueness implies that $\sigma ^{\ast}(f)=f$. Uniqueness also gives
injectivity of the morphism $A\rightarrow E_p$ and the above construction gives
surjectivity.
Now for every  point $p\in \pp ^1$ we obtain a complex structure $J_p$ on $A$
by pulling back the complex structure from $E_p$. The fact that $e\mapsto \sigma
(e)$ is antilinear means that $J_{\sigma p}=-J_p$.

Let $I=J_0$, $J=J_1$ and $K=J_{i}$.  We
claim that these provide a quaternionic triple and that for any other point $q$
the complex structure $J_q$ is that obtained from $(I,J,K)$ by the twistor space
construction described above. To prove this claim let $L$ be the
standard rank two pure twistor structure of weight $1$ with antipodal
involution,
corresponding to the quaternionic vector space ${\bf H}$.  The bundle
$\underline{Hom} (L, E)=L^{\ast}\otimes E$ is  pure of weight $0$ and has an
antipodal involution $\sigma$, thus it is the same as a vector space with real
structure. The real sections are the morphisms $L\rightarrow E$ compatible with
the antipodal involutions; let $Hom  (L,E)^{\sigma }$ denote this space of real
sections. It has real dimension equal to the complex dimension of $Hom (L,E)$
which in turn is the rank of $L^{\ast}\otimes E$ or $2r(E)$. If we fix a section
$u$ of $L$ which is preserved by $\sigma$ then evaluation at $u$ gives a
morphism
$$ \epsilon : Hom  (L,E)^{\sigma } \rightarrow A. $$
We claim that this morphism is injective. In fact if a $\sigma$-invariant map
$f:L\rightarrow E$ sends $u$ to zero then the bundle $ker(f)\subset L$ is
nontrivial, but it is a $\sigma$-invariant bundle, pure of weight $1$
so by the above discussion it corresponds to a nonzero subspace of $H^0(\pp ^1,
L)^{\sigma}={\bf H}$ which is preserved by all of the complex structures---hence
it must be all of ${\bf H}$ so $f=0$.

Now injectivity of the map $\epsilon$ implies surjectivity since the real
dimensions of the two sides coincide. Thus any section $a\in A$ is in the image
of a morphism $L\rightarrow E$. The commutation conditions for the various
complex structures on $L$ imply the same conditions for those complex structures
of $E$ as applied to $a$. This proves our claim that $(I,J,K)$ form a
quaternionic triple. The same proof shows that the other $J_q$ are as they are
supposed to be. It is obvious that $(E, \sigma )$ comes from this quaternionic
structure by the twistor construction.

There is a canonical antipodal and circular real twistor structure of
weight two which corresponds to the Tate Hodge structure.  We call it the ``Tate
twistor'' (!), denoted $T(1)$.  Tensorisation and dual leads to Tate twistors
$T(n)$ for all $n$.  By tensoring with these, any pure antipodal or
circular real twistor structure is equivalent to one of weight $0$ or
weight $1$,
so the above discussions apply: we obtain either a real vector space or a
quaternionic one.

Note, in passing, the fact that a complex vector space with quaternionic
structure must have even dimension. This carries over the Hodge-theory fact
that
a real Hodge structure of odd weight must have even dimension (on the other
hand, there is no such restriction for even weight).

The line bundle $\Oo (1)$ pure of weight $1$ (not to be confused with
$T(1)$ which is of weight two) has a natural circular real structure.  On the
dual $\Oo (-1)$ this is described as follows: over a point $P=[a:b] \in \pp ^1$
the fiber $\Oo (-1)_P$ is just the line of $(x,y)$ proportional to $(a,b)$.
Thus the total space of the bundle $\Oo (-1)$ is just ${\bf A}^2$ (blown up at
the origin).  We define the antilinear map
$$
\tau : {\bf A}^2 \rightarrow {\bf A}^2
$$
by
$$
\tau (x,y) = (\overline{y}, \overline{x}).
$$
This extends to the blow-up at the origin as an antilinear involution covering
the involution $\tau _{\pp ^1}$, thus giving a circular real structure on $\Oo
(-1)$. The dual is a circular real structure on $\Oo (1)$.  By tensoring
with a power of one of these, any pure twistor structure with circular real
structure becomes the same as one of weight zero, so it is just a real vector
space. This equivalence is compatible with taking the fiber $E_1$ over $1$
(which has a real structure since $1$ is a fixed point of $\tau _{\pp ^1}$).
Thus a circular real structure of any weight is just a real structure on the
underlying vector space.

\subnumero{Real structures in the Hodge case}

Let $H$ be the group of conformal automorphisms of $\pp  ^1$ which preserve
the set $\{ 0, \infty\}$ and which act on this set trivially when the
orientation is preserved, nontrivially when the orientation is changed.  The
connected component of $H$ is just $\Gm$ (the group of holomorphic automorphisms
fixing $0$ and $\infty$).  There are exactly two components, and the other
component is equal to $\Gm \cdot \sigma_{\pp ^1}$ or equally well $\Gm \cdot
\tau _{\pp ^1}$. In particular $H$ may be expressed as the group of
automorphisms generated by $\Gm$ and either one of $\sigma _{\pp ^1}$ or $\tau
_{\pp ^1}$.

\begin{proposition}
An $H$-equivariant mixed twistor structure $(E,W)$ is the same thing as an
$\rr$-mixed Hodge structure.
\end{proposition}
{\em Proof:}
Left to the reader.
\eop

This proposition says that mixed twistor structures which are both
$\Gm$-equivariant and have either an antipodal or circular real structure, are
real mixed Hodge structures. Since the groups generated by $\Gm$ and either of
$\sigma _{\pp ^1}$ or $\tau _{\pp ^1}$ are the same, we don't see the difference
between antipodal and circular real structures in the $\Gm$-equivariant (Hodge)
situation.

In the present note we shall not consider the problem of real twistor structures
any further (except one small place at the end when we discuss twistor spaces
for hypercomplex structures). The reader may imagine how to incorporate real
structures into all of the various statements we shall make: essentially all
constructions are equivariant for $\sigma$ or $\tau$ if these involutions are
given for the input data.

\numero{Variations of mixed twistor structure}
The problem of dealing with $\Cc ^{\infty}$ families of twistor structures
can be solved in several different ways.  The first would be to note that the
category of twistor structures has generators $\Oo _{\pp^1}(w)$ and that the
morphism spaces $Hom (\Oo _{\pp ^1}(w), \Oo _{\pp ^1}(w'))= \cc
[\lambda , \mu ]_{w'-w}$ are easily described as the spaces of homogeneous
polynomials (similarly for the extension groups).  One can then ``tensor'' this
category with any ring or sheaf of rings, in particular with $\Cc ^{\infty}_X$
for a manifold $X$.  Another more concrete approach would be to say that a
$\Cc ^{\infty}$ bundle of twistor structures over a manifold $X$ is a vector
bundle $E$ on $X\times \pp ^1$ provided with a complex structure operator
$\delbar _{\pp ^1}$ in the $\pp ^1$-direction.  For any $x\in X$ the restriction
$E|_{\{ x\} \times \pp ^1}$ becomes a holomorphic  vector bundle on $\pp ^1$,
algebraic by GAGA.  We will sometimes make use of this interpretation. Finally,
our main interpretation is sheaf-theoretic: if $X$ is a
$\Cc ^{\infty}$ manifold, let $\Cc ^{\infty}_X \Oo _{\pp ^1}$ be the sheaf of
functions on $X\times \pp ^1$ (in the usual topology) which are $\Cc^{\infty}$
and which are holomorphic in the $\pp ^1$-direction (in other words the
functions annihilated by the operator $\delbar _{\pp ^1}$ refered to above).
A {\em $\Cc ^{\infty}$ family of twistor structures on $X$} is now simply a
locally free sheaf $\Ee$ of $\Cc ^{\infty}_X \Oo _{\pp ^1}$-modules on $X\times
\pp ^1$. We often refer to this simply as a ``bundle'', or sometimes a
``$\Cc ^{\infty}_X \Oo _{\pp ^1}$-module'', on $X\times
\pp ^1$. A {\em strict filtration} of such a bundle is a filtration of $\Ee$ by
subsheaves such that locally on $X\times \pp ^1$ the filtration comes from a
decomposition of $\Ee$ into a direct sum of locally free $\Cc ^{\infty}_X \Oo
_{\pp ^1}$-modules. If $W_{\cdot}$ is a strict filtration then the
$Gr ^W_i(\Ee )$ are again locally free $\Cc ^{\infty}_X \Oo
_{\pp ^1}$-modules.

If $V$ is a $\Cc ^{\infty}$ bundle on $X$ provided with filtrations $F$ and $F'$
then by doing the Rees bundle construction in a $\Cc ^{\infty}$ family we obtain
a locally free $\Cc ^{\infty}_X \Oo _{\pp ^1}$-module $\xi (V, F, F')$ on
$X\times \pp ^1$.

The above constructions are functorial for $\Cc ^{\infty}$ morphisms
$X'\rightarrow X$, notably for inclusions of points $x\hookrightarrow X$.
Thus the above definitions restrict to the definitions from the beginning, over
each point $x\in X$.

A {\em $\Cc ^{\infty}$ family of mixed twistor structures over $X$} is a
$\Cc ^{\infty}_X \Oo _{\pp ^1}$-module $\Ee$ with strict filtration $W_{\cdot}$
such that the $Gr ^W_i(\Ee )$ are pure of weight $i$ (which means that for each
$x\in X$ the corresponding holomorphic bundle on $\{ x\} \times \pp ^1$ is
pure of
weight $i$).

Suppose $X$ is a complex manifold.
The cotangent bundle $T^{\ast}_{\cc}(X)$, which we often confuse with its sheaf
$\Aa ^{1}_X$ of sections, has a Hodge structure of
weight $1$ and consequently a twistor structure pure of weight $1$. This twistor
structure is a
$\Cc ^{\infty}_X \Oo _{\pp ^1}$-module
which is isomorphic to $T^{\ast}_{\cc}(X)
\otimes _{\cc} \Oo _{\pp
^1}(1)$. To be precise, we have a decomposition
$$
\Aa ^{1}_X = \Aa ^{1,0}_X \oplus \Aa ^{0,1}_X,
$$
and the twistor structure which is $\xi (\Aa ^{1}_X, F, F')$ decomposes as
$$
\xi (\Aa ^{1}_X, F, F')= \xi (\Aa ^{1,0}_X, F, F')\oplus \xi (\Aa ^{0,1}_X, F,
F'),
$$
where now the filtrations on the pieces
$\Aa ^{1,0}_X$ and $\Aa ^{0,1}_X$ are trivial (but shifted differently).
These trivial filtrations induce isomorphisms
$$
\xi (\Aa ^{1,0}_X, F, F')\cong
\Aa ^{1,0}_X\otimes _{\cc} \Oo _{\pp^1}(1),\;\;\;
\xi (\Aa ^{0,1}_X, F, F')\cong
\Aa ^{0,1}_X\otimes _{\cc} \Oo _{\pp^1}(1)
$$
and their sum is the isomorphism
$$
\phi :\xi (\Aa ^{1}_X, F, F')\cong
\Aa ^{1}_X\otimes _{\cc} \Oo _{\pp^1}(1)
$$
in question.
There is a natural inclusion of bundles
$$
\iota : \Aa ^{1}_X\otimes _{\cc} \Oo _{\pp
^1}\hookrightarrow \xi (\Aa ^{1}_X, F, F')
$$
and this decomposes as inclusions on each of the pieces. If we fix the standard
sections $\lambda$ and $\mu$ of of $\Oo (1)$ which vanish
respectively at $0$ and $\infty$, we can write
$$
\phi \iota (\alpha ^{1,0} + \alpha ^{0,1})=
\lambda \alpha ^{1,0} + \mu \alpha ^{0,1}.
$$
The inclusion $\iota$ fixes notation for $\xi (\Aa ^{1}_X, F, F')$
and the above formula fixes notation for the isomorphism  $\phi$.

We denote for short $\xi (\Aa ^{1}_X, F, F')$ by $\xi \Aa ^1_X$. The same
construction applies to the whole exterior algebra $\Aa ^{\cdot}_X$
to give a graded $\Cc ^{\infty}_X\Oo _{\pp ^1}$-module $\xi\Aa ^{\cdot}_X$
on $X\times \pp ^1$. The $\xi\Aa ^i_X$ have twistor structures  of weight $i$.

There is a natural differential operator
$$
{\bf d}: \xi\Aa ^0_X \rightarrow \xi\Aa ^1_X
$$
which is a morphism of twistor structures.  It is just the composition of
the usual exterior derivative $d$ followed by $\iota$, so it can be written
$$
{\bf d} = \lambda \partial  + \mu \delbar
$$
where again $\lambda$ and $\mu$ are the sections of $\Oo _{\pp ^1}(1)$ which
vanish respectively at $0$ and $\infty$.  The differential extends to the
comples  $\xi\Aa ^i_X$ have twistor structures  of weight $i$ and we obtain a
complex $\xi\Aa ^{\cdot}_X$ with differential ${\bf d}$.

We generalize the definition of \cite{SteenbrinkZucker} to the twistor case.
A {\em variation of mixed twistor structure}
\footnote{This is the twistor analogue of the complex variations of
mixed Hodge structure; one could also talk about antipodal or circular real
variations of mixed twistor structure which would be analogues of $\rr$-VMHS,
this is left to the reader.}
is a $\Cc ^{\infty}$ family of mixed twistor structures  $(\Ee ,
W_{\cdot})$ on $X$ (i.e. a
$\Cc ^{\infty}_X\Oo _{\pp ^1}$-module $\Ee$ with strict filtration
$W_{\cdot}$ on $X\times \pp ^1$)
together with an operator
$$
D: \Ee \rightarrow \Ee \otimes
_{\Cc ^{\infty}_X\Oo _{\pp ^1}}\xi \Aa ^1_X
$$
respecting the weight filtration, such that
the Leibniz rule
$$
D(ae) = {\bf d}(a) e + aD(e)
$$
is satisfied, and such that $D^2=0$.
A {\em pure variation of twistor structure} is the same as above but with
the associated graded of the weight filtration concentrated in one degree.

Suppose given a variation of mixed twistor structure.
For any $\lambda \in {\bf A}^1$ we obtain an {\em underlying
$\lambda$-connection}
by taking the fiber over $\lambda \in \pp ^1$.
Note that the fiber of $\xi \Aa ^1_X$ over $\lambda \in \pp ^1$ is naturally
identified with $\Aa ^1_X$ and via this identification, ${\bf d}$ corresponds to
$\lambda \partial + \delbar$.   If $(\Ee , D)$ is a variation of twistor
structure then the bundle $\Ee _{\lambda}:= \Ee |_{X\times \{ \lambda \} }$
is a $\Cc ^{\infty}$ bundle on $X$ with operator $D_{\lambda}$ having symbol
$\lambda \partial + \delbar$ and square zero.  If we decompose according to
Hodge type
$$
D_{\lambda}= D^{1,0}_{\lambda} + D^{0,1}_{\lambda}
$$
then $D^{0,1}_{\lambda}$ provides an integrable holomorphic structure for $\Ee
_{\lambda}$ and $D^{1,0}_{\lambda}$ becomes a holomorphic $\lambda$-connection
(i.e. operator satisfying Leibniz for $\lambda d$ and having square zero---cf
\cite{SantaCruz}; this definition was made by Deligne in \cite{DeligneLetter}).

In particular for $\lambda = 0$ we obtain an
underlying Higgs bundle and for $\lambda = 1$ we obtain an underlying flat
bundle. The $\lambda$-connections come in a holomorphic family
indexed by ${\bf A}^1$, as will be explained in further detail below.

On a slightly different note, for every point $x\in X$ we obtain a mixed twistor
structure $(\Ee, W_{\cdot})_x$. This gives a map $X\rightarrow {\cal MTS}$
to the
moduli stack of mixed twistor structures; note however that it is only a
$\Cc^{\infty}$ map and not holomorphic. We call this the {\em classifying map}
for the variation of mixed twistor structure $(\Ee, W_{\cdot})$.

\subnumero{Polarizations}
We define a notion of {\em polarization} for pure variations of twistor
structures.  This will be used as a  characterization of those variations which
correspond to harmonic bundles, cf Lemma \ref{harmonic} below.  As we will not
use the notion of polarization anywhere else, the reader may prefer to skip
directly to Lemma \ref{harmonic} and take as the definition of polarizable
variation, one which corresponds to a harmonic bundle by the construction of
\ref{harmonic}. The {\em raison d'etre} of the definition we give below of
polarization is just to emphasize the analogy with variations of Hodge
structure.

Suppose $(\Ee , {\bf d})$ is a pure variation of twistor structure of weight
$w$ on $X$.
We obtain the locally free sheaf of  $\Cc ^{\infty}_{X}\Oo _{\pp
^1}$-modules $\sigma ^{\ast}(\Ee )$ on $X\times \pp ^1$ as
follows:
if $U\subset X\times \pp ^1$ is an open set then $(1\times \sigma )(U)$
is an open subset of $X\times \pp ^1$ (and all open subsets are
obtained this way).  We set
$$
\sigma ^{\ast}(\Ee )((1\times \sigma )(U)):=\Ee (U).
$$
We give this a structure of $\Cc ^{\infty}_{X}\Oo _{\pp
^1}$-module as follows: if
$e\in \Ee (U)$ and
$$
a\in \Cc ^{\infty}_{X}\Oo _{\pp
^1}((1\times \sigma )(U))
$$
then $a$ times $e$ is defined to be equal to
$$
\overline{(1\times \sigma )^{\ast}(a)} e\in \Ee (U).
$$
The complex conjugate is required in order to obtain a section
$\overline{(1\times \sigma )^{\ast}(a)}$ which is holomorphic in the $\pp
^1$-direction. This has the effect that the functor $\Ee \mapsto
\sigma ^{\ast}(\Ee )$ is a $\cc$-antilinear functor.

We give $\sigma  ^{\ast}(\Ee )$  a structure of variation of
twistor structure on $\overline{X}$. For this, we define an operator
$$
\sigma ^{\ast}(D): \sigma  ^{\ast}(\Ee )
\rightarrow \sigma  ^{\ast}(\Ee )
\otimes
_{\Cc ^{\infty}_{X}\Oo _{\pp ^1}}\xi\Aa ^1_X
$$
by using a morphism
$$
\sigma ^{\ast}(\xi\Aa ^1_X)\rightarrow \xi\Aa ^1_X.
$$
This morphism is defined as follows. Sections of $\xi\Aa ^1_X$ may be written in
the form $\alpha ' \lambda + \alpha '' \mu$ where $\alpha '$ and $\alpha ''$
are functions on $X\times \pp ^1$ taking values respectively in the $1,0$ and
$0,1$ forms on $X$, and holomorphically varying in the $\pp ^1$-direction.
Furthermore $\alpha '$ is allowed to have one pole at the zero of $\lambda$
(i.e. along $X\times \{ 0\}$) and $\alpha '' $ is allowed to have one pole
at the
zero of $\mu$ (i.e. along $X\times \{ \infty \}$).   Noting that
$\sigma ^{\ast}(\xi\Aa ^1_X)(U) = \xi \Aa ^1_X ((1\times \sigma )(U))$,
we define the morphism
$$
\psi :\xi \Aa ^1_X ((1\times \sigma )(U))\rightarrow \xi \Aa ^1_X(U)
$$
by
$$
\psi (\alpha ' \lambda + \alpha '' \mu ):= \overline{\sigma ^{\ast}(\alpha ''
)}\lambda + \overline{\sigma ^{\ast}(\alpha '
)}\mu .
$$
This expression still satisfies the conditios described above concerning the
allowable poles of the coefficients, since $\sigma$ interchanges $0$ and
$\infty$. Using this morphism we obtain the operator $\sigma ^{\ast}(D)$,
making $(\sigma ^{\ast}(\Ee ), \sigma ^{\ast}(D))$ into a variation of twistor
structure.

We note in passing that there isn't really any other way to define the antipodal
conjugate $(\sigma ^{\ast}(\Ee ), \sigma ^{\ast}(D))$; this is a contrast with
the situation in Hodge theory where one doesn't really see any inner logic
requiring the various changes of signs. (In a certain sense we have transfered
the question to the simple fact of choosing to work with the antipodal
involution, which turns out to contain all of the necessary sign changes.)

A {\em polarization} of $(\Ee , {\bf d})$ is
a bilinear pairing
$$
P: \Ee \otimes \sigma ^{\ast} (\Ee ) \rightarrow T(w)
$$
which is a morphism of variation of twistor structures, and which is {\em
positive hermitian}. This last notion is defined as follows: the form $P$
induces
a morphism of trivial bundles
$$
\Ee (-w) \otimes \sigma ^{\ast}(\Ee (-w)) \rightarrow T(w) \otimes \Oo (-w)
\otimes \Oo (-w) \cong \Oo ,
$$
hence a morphism on the corresponding vector spaces; but
$H^0(\sigma ^{\ast}(\Ee (-w)) )\cong H^0(\Ee (-w))$ by an antilinear
isomorphism,
and we say that $P$ is {\em hermitian} (resp. {\em positive hermitian}) if the
resulting antilinear form on the vector space $H^0(\Ee (-w))$ is hermitian
(resp. positive hermitian). Note that the hermitian condition can be expressed
in sheaf-theoretic terms without refering to $H^0$ but I don't see a nice way to
do this for the positivity condition.

We say that the variation $(\Ee , {\bf d})$ is
{\em polarizable} if there exists a polarization.

We say that a variation of mixed twistor structure $(\Ee , W_{\cdot}, {\bf d})$
is {\em graded-polarizable} if the associated-graded pieces $Gr ^W_i(\Ee )$ are
polarizable as pure variations of twistor structure.

\begin{lemma}
\label{harmonic}
If $V$ is a flat bundle with pluriharmonic  metric then $V$ underlies a
polarizable pure variation of twistor structure. In particular, if $X$ is a
compact K\"{a}hler manifold then  any irreducible representation of $\pi _1(X)$
corresponds to a flat bundle underlying a pure variation of twistor structure.
This variation is unique up to change of weight (by tensoring with a one
dimensional twistor structure). Any polarizable variation of pure twistor
structure is a direct sum of ones obtained from irreducible representations
by this construction.
\end{lemma}
{\em Proof:}
Suppose $\Ee$ is a pure variation of twistor structure of weight zero (the other
weights are treated by tensoring with $\Oo _{\pp ^1}(w)$).  Let $p: X\times \pp
^1\rightarrow X$ denote the first projection. Then $E:=p_{\ast}(\Ee )$
is a locally free $\Cc^{\infty}_X$-module on $X$ of rank equal to the rank of
$\Ee$, and we have
$$
\Ee = p^{-1}(E)\otimes _{p^{-1}\Cc^{\infty}_X}
\Cc^{\infty}_X\Oo _{\pp ^1}.
$$
Similarly
$p_{\ast} \xi \Aa ^0_X = \Aa ^0_X$. However,
$$
p_{\ast} \xi \Aa ^1_X \cong \Aa ^1_X \otimes _{\cc} \cc \langle \lambda , \mu
\rangle
$$
is a bundle of rank twice that of $\Aa ^1_X$,
since
$\xi \Aa ^1_X \cong \Aa ^1_X \otimes \Oo _{\pp ^1} (1)$ as described previously.
We have
$$
p_{\ast}({\bf d}) = \lambda \partial + \mu \delbar ,
$$
and $p_{\ast}(D)$ is an operator
$$
p_{\ast}(D): E\rightarrow E \otimes _{\Cc^{\infty}_X}
\Aa ^1_X \otimes _{\cc} \cc \langle \lambda , \mu
\rangle
$$
whose symbol is $\lambda \partial + \mu \delbar$.  We can write
$$
p_{\ast}(D) = \lambda D' + \mu D''
$$
where $D'$ and $D''$ are operators from $E$ to $E\otimes _{\Cc^{\infty}_X}
\Aa ^1_X$. These operators are uniquely determined by the above equation
and furthermore it follows that $D'$ has symbol $\partial$ and $D''$ has symbol
$\delbar$. The condition $D^2=0$ implies the integrability conditions
$(D')^2=0$,
$(D'' ) ^2=0$ and $D'D'' + D'' D' = 0$.

Conversely suppose $(E,D', D'')$ is a bundle with operators $D'$ and $D''$
with symbols $\partial$ and $\delbar$ respectively and
satisfying the above integrability conditions. Then setting $$ \Ee :=
p^{-1}(E)\otimes _{p^{-1}\Cc^{\infty}_X} \Cc^{\infty}_X\Oo _{\pp ^1},
$$
we get a $\Cc ^{\infty}$ family of pure twistor structures of weight zero.
Putting
$D:= \lambda D' + \mu D''$ considered as an operator from
$\Ee$ to $\Ee \otimes _{\Cc^{\infty}_X\Oo _{\pp ^1}} \xi \Aa ^1_X$
via the previously-mentionned isomorphism $\xi \Aa ^1_X\cong \Aa ^1_X
\otimes _{\cc } \Oo _{\pp ^1}(1)$ (and via consideration of $\lambda $ and $\mu$
as sections of $\Oo _{\pp ^1}(1)$), we obtain a variation of pure twistor
structure $(\Ee , D)$.

These constructions establish a one to one correspondence between pure
variations of twistor structure $(\Ee , D)$ of weight zero, and triples $(E,D',
D'')$ as in the definition of harmonic bundle \cite{HBLS}.  We claim that
the variation $(\Ee , D)$ is polarizable if and only if the operators $D'$ and
$D''$ are related by a metric (which is thus a harmonic metric) $K$ on $E$
according to the definitions in \cite{HBLS}.

The first thing to note is that under the above correspondence (assuming $\Ee$
is of weight zero) we have that $p_{\ast}(\sigma ^{\ast}\Ee )= \overline{E}$ is
the complex conjugate $\Cc ^{\infty}$ bundle. We have to calculate
$p_{\ast}(\sigma ^{\ast} D)$.  Note that we can decompose $D'$ and $D''$
according to Hodge type of forms and write
$$
p_{\ast}(D) = \lambda D'_{1,0} + \lambda D' _{0,1} +
\mu D''_{1,0} + \mu D'' _{0,1} .
$$
If we write this in a form ready to apply the morphism $\psi$ used above to
define $\sigma ^{\ast}(D)$ it becomes
$$
p_{\ast}(D) = \lambda (D'_{1,0} + \frac{\mu }{\lambda}D''_{1,0} )
+ \mu (
\frac{\lambda}{\mu} D' _{0,1} +
 D'' _{0,1} ).
$$
Now $p_{\ast}(\sigma ^{\ast} D)$ is obtained by applying $\sigma ^{\ast}$ to the
$\Ee$-coefficients, and by applying the operation $\psi$ to the
form coefficients and $\lambda $ and $\mu$. We write this as
$$
p_{\ast}(\sigma ^{\ast} D)= \lambda \cdot \overline{\sigma ^{\ast}
( \frac{\lambda}{\mu} D' _{0,1} +
 D'' _{0,1} )} +
\mu \cdot \overline{\sigma ^{\ast}  (D'_{1,0} + \frac{\mu
}{\lambda}D''_{1,0} )}.
$$
The pieces $D'_{1,0}$ and so on are invariant under $\sigma ^{\ast}$ because
they are constant in the $\pp ^1$-direction. However, note that
$$
\overline{\sigma ^{\ast}(\frac{\lambda}{\mu})}
=-\frac{\mu}{\lambda}
$$
since the antipodal involution is written $t\mapsto -\overline{t}^{-1}$
in terms of the coordinate $t$ on ${\bf A}^1$. Thus we have
$$
p_{\ast}(\sigma ^{\ast} D)= \lambda \cdot
( -\frac{\mu}{\lambda} \overline{D}' _{0,1} +
 \overline{D}'' _{0,1} ) +
\mu \cdot  (\overline{D}'_{1,0} - \frac{\lambda }{\mu}\overline{D}''_{1,0} )
$$
$$
= \lambda (\overline{D}'' _{0,1}- \overline{D}''_{1,0}) + \mu
(\overline{D}'_{1,0}-\overline{D}' _{0,1}).
$$
In other words the decomposition
$$
p_{\ast}(\sigma ^{\ast} D)= \lambda (\sigma ^{\ast}D)' + \mu (\sigma
^{\ast} D)''
$$
is given by
$$
(\sigma ^{\ast}D)' = \overline{D}'' _{0,1}- \overline{D}''_{1,0}
$$
and
$$
(\sigma ^{\ast}D)'' = \overline{D}'_{1,0}-\overline{D}' _{0,1}.
$$
Thus the triple  associated to $(\sigma ^{\ast}\Ee , \sigma ^{\ast}D)$
is $(\overline{E}, \overline{D}'' _{0,1}- \overline{D}''_{1,0},
\overline{D}'' _{0,1}- \overline{D}''_{1,0})$.

A polarization $P$ of $\Ee$ corresponds to a  morphism
$$
K:= p_{\ast}(P): E\otimes \overline{E} \rightarrow \cc
$$
which, by hypothesis, is a positive definite hermitian form on $E$. The
morphism $K$ intertwines the pair of operators $(D', D'')$ on $E$ with the pair
of operators $(overline{D}'' _{0,1}- \overline{D}''_{1,0},
\overline{D}'' _{0,1}- \overline{D}''_{1,0})$ on $\overline{E}$---which exactly
says that $K$ is a harmonic metric for the triple $(E,D', D'')$. Conversely by
following the above formulas in the other direction, a harmonic metric leads to
a polarization which is a morphism of mixed twistor structures.
This completes the proof of the claimed correspondence between polarizations and
harmonic metrics.

This claim implies the lemma, modulo the remark that we can pass from
variations $(\Ee , D)$ pure of weight $w$ to variations pure of weight zero (and
back again) by tensorization by the constant rank one twistor
structure $\Oo _{\pp
^1}(-w)$ (resp.   $\Oo _{\pp
^1}(w)$). Note that these twistor structures admit polarizations, so tensoring
with them preserves the polarization condition.
\eop

{\em Remark:}  This lemma is the fundamental reason for the utility of the
notion of twistor structure: it applies to any irreducible representation of the
fundamental group, and not just certain ones (the variations of Hodge
structure).

\subnumero{Description in terms of $t$-connections on $X$ and
$\overline{X}$}

On $X\times {\bf A}^1 \subset X \times \pp ^1$ we have an isomorphism
$$
\Phi : \Aa ^i_X |_{X\times {\bf A}^1 } \cong A^i (X\times {\bf A}^1/ {\bf A}^1).
$$
Via this isomorphism we have
$$
\Phi \circ {\bf d} = \delbar + t \partial
$$
where $t$ is the coordinate on $\pp ^1$.
Using this isomorphism we can identify an operator $D$ with an operator
$$
D_0 : \Ee \rightarrow \Ee \otimes A^1(X\times {\bf A}^1/ {\bf A}^1)
$$
satisfying the Leibniz rule
$$
D_0 (ae) = aD_0(e) + (\delbar a + t \partial a) e.
$$
Decompose according to type as $D_0 = D_0^{1,0}+ D_0^{0,1}$. Furthermore we can
add the operator $\partial $ in the ${\bf A} ^1$-direction to obtain an operator
$$
D_0^{0,1} + \delbar : \Ee \rightarrow A^{0,1}(X\times {\bf A}^1, \Ee ).
$$
This operator satisfies Leibniz' rule for the symbol $\partial_{X\times {\bf
A}^1}$ and it is integrable (since $D_0^2= 0$, $\delbar ^2=0$ and
$[\delbar , D_0]=0$). Thus this operator defines a holomorphic structure for the
bundle $\Ee|_{X\times {\bf A}^1}$.  Denote this holomorphic bundle---or more
precisely its sheaf of holomorphic sections---by $\Ff$. The remaining part
$D^{1,0}_0$ commutes with the holomorphic structure and satisfies Leibniz' rule
for the symbol $t\partial$.  In other words this operator leads to a
morphism of
holomorphic sheaves
$$
\nabla : \Ff \rightarrow \Ff \otimes _{\Oo _{X\times {\bf A}^1}} \Omega
^1_{X\times {\bf A}^1/{\bf A}^1}
$$
satisfying the Leibniz rule $\nabla (af)= a\nabla (f) + d(a)\nabla (f)$
for  sections $a$ of $\Oo _{X\times {\bf A}^1}$ and $f$ of $\Ff$.
Furthermore $\nabla ^2=0$. Thus $(\Ff , \nabla )$ is a {\em $t$-connection on
$X\times {\bf A}^1/{\bf A}^1$}.

We obtain by symmetry a similar description on the other
standard affine open set
${\bf A}^1\subset
\pp ^1$ whose coordinate is $t^{-1}$, but with $X$ replaced by $\overline{X}$.
We are forced to replace $X$ by $\overline{X}$ because it is the $1,0$ part of
the connection on $X$ whose symbol doesn't degenerate. This yields a
$t^{-1}$-connection $(\Ff ', \nabla ')$ on  $\overline{X}\times {\bf A}^1/{\bf
A}^1$.

In order to relate these two descriptions, in other words to understand the
glueing between the objects $(\Ff , \nabla )$ and $(\Ff ', \nabla ')$
we need a different trivialization of $\Aa ^i_X$ over $X\times \Gm \subset X
\times \pp ^1$,
$$
\Psi : \Aa ^1_X |_{X\times \Gm} \cong A^1(X\times \Gm / \Gm )
$$
such that
$$
\Psi \circ {\bf d} = d = \delbar + \partial .
$$
Via this trivialization, $\Ee |_{X\times \Gm}$ becomes a vector bundle with a
flat connection relative to the base $X\times \Gm$ and a commuting operator
$\delbar$ in the horizontal direction. The sheaf of sections annihilated by
the connection and the $\delbar$-operator is a locally free sheaf of
$p_1^{-1}(\Oo _{\Gm})$-modules on $X^{\rm top}\times \Gm$, in other words it is
a holomorphic family of local systems which we denote $\Ll = \{ L_t\} _{t\in
\Gm}$ on $X^{\rm top}$.  (If we choose a basepoint $x\in X$ then this is
the same
as a holomorphic family of representations of $\pi _1(X,x)$ modulo holomorphic
changes of basis.)

The relationship between $(\Ff , \nabla )|_{X\times \Gm}$ and $\Ll$ is that
for any $t\in \Gm$, $t^{-1}\nabla$ is a holomorphic flat connection and $L_t$ is
its local system of flat sections.
By symmetry a similar interpretation holds for $(\Ff ', \nabla ')$.

\numero{Cohomology of smooth compact K\"{a}hler manifolds with VMTS
coefficients}

Classically, the cohomology of a smooth compact K\"{a}hler manifold with
constant
coefficients was the first object to be provided with a pure Hodge structure.
Applying our meta-theorem to this statement doesn't yield anything other than
trading in the Hodge structure for a twistor structure, since the notion of
Hodge structure doesn't appear in the hypothesis.  One of the main starting
points for modern Hodge theory was Deligne's observation that the K\"ahler
identities also work with coefficients in a variation of Hodge structure
\cite{DeligneLetSerre}. This
leads to the theorem that the $k$-th cohomology of a smooth compact K\"ahler
manifold with coefficients in a pure variation of Hodge structure of weight $n$,
carries a pure Hodge structure of weight $n+k$.  This was generalized slightly
in \cite{SteenbrinkZucker} to the case of coefficients in a variation of mixed
Hodge structure, giving a mixed Hodge structure on the cohomology. We generalize
to twistor structures in this section.

\begin{theorem}
\label{coho1}
{\rm (\cite{DeligneLetSerre}, \cite{SteenbrinkZucker})}
Suppose $(\Ee , W_{\cdot} , {\bf d})$ is a $Gr$-polarizable variation of mixed
twistor structure on a compact K\"ahler manifold $X$. Then $H^k(X, \Ee _1)$ (the
cohomology with coefficients in the underlying flat bundle) carries a natural
mixed twistor structure. If the variation is pure of weight $n$ then the $k$-th
cohomology is pure of weight $n+k$.
\end{theorem}

\subnumero{Pure coefficients}
First we treat the case of pure coefficients.
Suppose that $X$ is a smooth compact K\"{a}hler variety and
$(\Ee ,{\bf d})$ is a polarizable variation of pure twistor structure of weight
$n$ on $X$.   Let $L_1$ denote the underlying flat bundle (it is the fiber over
$1$ of the associated family of flat bundles $\{ L_t\}$).  We will construct a
pure  twistor structure of weight $n+k$ on the $k$-th cohomology of $X$ with
coefficients in $L_1$.

In this case note that there is a harmonic bundle $(E, D', D'')$ and
$$
\Ee= p_1^{\ast} E \otimes p_2^{\ast}\Oo _{\pp ^1}(n),
$$
with ${\bf d}= \lambda D' +
\mu D'' $. Here the pullbacks $p_1^{\ast}$ and $p_2^{\ast}$ take $\Cc
^{\infty}_X$-modules (resp. $\Oo _{\pp ^1}$-modules) to
$\cc ^{\infty}_X\Oo _{\pp ^1}$-modules on $X\times \pp ^1$.

Let $\xi \Aa ^{\cdot}(\Ee )$ be the twistor complex of forms with
coefficients in
$\Ee$,
a complex of $\Cc ^{\infty}_X\Oo _{\pp ^1}$-modules on $X\times \pp ^1$ with
operator ${\bf d}$. We have
$$
\xi \Aa ^i(\Ee ) \cong p_1^{\ast}(A^i(E)) \otimes p_2^{\ast}\Oo _{\pp ^1}(n+i).
$$
Let $H^i(E) \subset A^i(E)$ denote the subspace of harmonic $i$-forms
with coefficients in $E$. We have for any $(a,b)\neq (0,0)\in \cc ^2$
isomorphisms
$$
H^i(E) \stackrel{\cong}{\rightarrow} H^i(A^{\cdot}(E), aD' + bD'').
$$
This is shown in \cite{HBLS} for $(a,b) = (0,1)$ and $(1,1)$, but the same proof
(which is just a direct generalization of the standard utilization of the
K\"ahler identities) works in general---one defines the laplacian for
$aD' + bD''$ and shows that it is proportional to the laplacian for $D''$. We
get back to Deligne's original idea about extending the K\"ahler identities.
In
particular, the above morphism induces an isomorphism
$$
H^i(E) \otimes   \Oo _{\pp ^1}(n+i)
\stackrel{\cong}{\rightarrow} M:= R^ip_{2,\ast}(\xi \Aa ^{\cdot}(\Ee ),
{\bf d}).
$$
This proves that the higher direct image is a bundle $M$ which is  pure of
weight
$n+i$.

We  now give a second, analytic construction of $M$, which will be useful later
on. However, we refer to the above differential-geometric point of view to prove
purity of the weight quotients of our structure.

Let $\Ff$, $\Ll$ and $\Ff '$
be the triple associated to $\Ee$ consisting of a $t$-connection $(\Ff ,
\nabla )$
on $X\times {\bf A}^1$, a family of local systems $\Ll = \{ L_t\}$ and a
$t^{-1}$-connection $(\Ff ', \nabla ')$ on the other $\overline{X}\times {\bf
A}^1$. Recall that  $L_t$ is the flat bundle over $X^{\rm top}$ associated
to the
holomorphic flat connection  $(\Ff _t, t^{-1} \nabla _t)$ on $X$, and it is also
the flat bundle associated to the  connection
$(\Ff' _{t^{-1}}, t \nabla _{t^{-1}})$ on $\overline{X}$.

Let $\xi \Omega ^{\cdot}_{X}(\Ff )$
denote the Rees bundle complex of the complex of
relative differentials with coefficients in $\Ff$, with differential given
by the $t$-connection $\nabla$.  More precisely we have
$$
\xi \Omega ^{\cdot}_{X}(\Ff )= \xi (\Omega ^{\cdot}_X, F) \otimes _{\Oo
_{X\times \pp ^1}}\Ff ,
$$
where the filtration $F$ is the usual Hodge filtration (i.e. the ``stupid
filtration'').

{\em Remark:} there shouldn't be too much confusion between our use of the
symbol $\xi$ here for the Rees bundle construction with one filtration, and its
use elsewhere for the Rees bundle construction with two filtrations: it depends
on whether we are working over the affine line (as is presently the case) or
over $\pp ^1$ (as was the case in the differential-geometric construction
above).

For $t\neq 0$ the hypercohomology of
$\xi \Omega ^{\cdot}_{X}(\Ff )|_{X\times \{ t\} }$
calculates the cohomology of the local system $L_t$.
Thus over
$\Gm$ we have
$$
{\bf R}^kp_{2,\ast}(\xi \Omega ^{\cdot}_{X}(\Ff )
|_{X\times \Gm} \cong
{\bf R}^kp_{2,\ast}\Ll .
$$

Let $\xi \Omega ^{\cdot}_{\overline{X}}(\Ff ')$
denote the Rees bundle complex of the complex of
relative differentials with coefficients in $\Ff'$, with differential given
by the $t^{-1}$-connection $\nabla '$. Note that this is now taking place
over the affine line neighborhood of $\infty$. Again, over $\Gm$ we have
$$
{\bf R}^kp_{2,\ast}(\xi \Omega ^{\cdot}_{\overline{X}} (\Ff ')
|_{\overline{X}\times \Gm} \cong
{\bf R}^kp_{2,\ast}\Ll .
$$
We can use these two isomorphisms to glue together
${\bf R}^kp_{2,\ast}(\xi \Omega ^{\cdot}_{X}(\Ff )$
and ${\bf R}^kp_{2,\ast}(\xi \Omega ^{\cdot}_{\overline{X}} (\Ff ')$ to obtain a
bundle $M$ over $\pp ^1$.   This bundle is the same as that constructed
previously; in particular it is pure of weight $n+k$.

\subnumero{Mixed coefficients}
For this section suppose that $X$ is a smooth compact K\"{a}hler variety and
$(E,W,{\bf d})$ is a graded-polarizable variation of mixed twistor structure on
$X$.   Let $L_1$ denote the underlying flat bundle (it is the fiber over $1$ of
the associated family of flat bundles $\{ L_t\}$).  We will construct a mixed
twistor structure  on the $k$-th cohomology of $X$ with coefficients in $L_1$
(fix $k$ for the rest of this section).

Let $\xi \Aa ^i_X(E)$ be the twistor complex of forms with coefficients in
$E$. It
is a complex of locally free $\Cc ^{\infty}_X\Oo _{\pp ^1}$-modules on $X\times
\pp ^1$ with operator ${\bf d}$.

Define the {\em pre-weight filtration} as
$$
W^{\rm pre}_n \xi \Aa^i_X(E):= \xi \Aa^i_X(W_nE).
$$
Let $M$ denote the $k$-th cohomology sheaf on $\pp ^1$ of the filtered
complex
$$
\Mm ^{\cdot} := p_{2,\ast} (\xi \Aa^i_X(E), W^{\rm pre}, {\bf d}).
$$
Note that $M$ has a filtration which we denote $W^{\rm pre}M$ and we define
$$
W_mM := W^{\rm pre}_{m-k}M .
$$
We claim that $(M,W)$ is  a mixed
twistor structure. Given the previous result in the pure case,
this is essentially
the same as Lemma \ref{degen} below but we give the direct argument here as a
warmup. The spectral sequence for the  cohomology of the filtered complex $\Mm
^{\cdot}$ starts with  $$
E^{p,q}_1 = H^{p+q} (Gr ^{W^{\rm pre}}_{-p}(\xi \Aa^{\cdot}_X(E)) )
$$
$$
= H^{p+q} (\xi \Aa ^{\cdot}_X (Gr ^W_{-p}(E))),
$$
which is a pure twistor structure of weight $p+q-p = q$ by the above result for
the pure polarized variations of twistor structure $Gr ^W_{-p}(E)$.  The
differential $d_1:E^{p,q}_1 \rightarrow E^{p+1,q}_1$ is a morphism of twistor
structures of the same weight, so the cohomology of the differential $d_1$ are
pure twistor structures $E^{p,q}_2$ of weight $q$. Now the differential $d_2:
E^{p,q}_2 \rightarrow E^{p+2, q-1}_2$ must vanish since it is a morphism from a
semistable bundle of slope $q$ to a semistable bundle of slope $q-1$. Thus
$E^{p,q}_3=E^{p,q}_2$ and arguing by induction, $E^{p,q}_2=E^{p,q}_{\infty}$ are
pure twistor structures of weight $q$. But
$$
E^{p,q}_{\infty} = Gr ^{W^{\rm pre}}_{-p}H^{p+q}(\xi \Aa^{\cdot}_X(E)),
$$
in particular
$$
Gr ^W_m (M)= Gr ^{W^{\rm pre}}_{m-k}H^{k}(\xi \Aa^{\cdot}_X(E))
= E^{k-m, m}_{\infty}.
$$
This shows that $Gr ^W_m (M)$ is pure of weight $m$ as desired.

One can do an analytic construction of $M$ parallel to the analytic construction
in the pure case. This would essentially be repeating what we will say below
concerning patching, so for now it is left to the reader.

We have now proved Theorem \ref{coho1}.
\eop

\numero{Cohomology of open and singular varieties}
We generalize \cite{Hodge2}.  For now we stay with the easy situation of
coefficients which extend across a compactification, and don't attack the
problem
of coefficients which are a local system only defined on the open
variety---since even in the case of variations of Hodge structure, this is a
much harder problem.

\begin{theorem}
\label{TwistorII}
{\rm (cf \cite{Hodge2}, 3.2.5)}
Suppose $X$ is a compact K\"{a}hler manifold with a divisor $D\subset X$ with
normal crossings. Let $U= X-D$. Suppose $(\Ee , W_{\cdot}, D)$ is a variation
of mixed twistor structure on $X$ with underlying flat bundle $V$.   Then the
cohomology $H^i(U, V|_U)$ carries a natural mixed twistor structure, which can
be described as the higher direct image of the de Rham complex of
$(\Ee , W_{\cdot}, D)|_U$ via the projection $U\times \pp ^1 \rightarrow
\pp ^1$.
This mixed twistor structure is functorial in $U$ (independent of the
compactification).
\end{theorem}

We adopt the following definition: if $X$ is any complex analytic space, then a
{\em variation of mixed twistor structure on $X$} is a functorial assignment,
for every smooth complex manifold $M$ mapping to $X$, of a variation of mixed
twistor structure on $M$.  To be precise about this, it means that for
every morphism $f: M\rightarrow X$ from  a smooth complex manifold we have
a VMTS
${\bf E}(f) = (\Ee (f), W_{\cdot}(f), {\bf d}(f))$ on $M$; and whenever $a:
M\rightarrow M'$ and $f':M'\rightarrow X$ with $f=f'a$ then we have an
isomorphism $\epsilon (a): a^{\ast}({\bf E}(f'))\cong {\bf E}(f)$, such
that the cocycle
condition holds: if  $M\stackrel{a}{\rightarrow}
M'\stackrel{b}{\rightarrow} M''$
and $f'': M'' \rightarrow X$ with $f'=f''b$ and $f=f''ba$ then the composition
$$
a^{\ast} b^{\ast}({\bf E}(f'')) \stackrel{a^{\ast}\epsilon (b)}{\rightarrow}
a^{\ast}({\bf E}(f'))\stackrel{\epsilon (a)}{\rightarrow} {\bf E}(f)
$$
is equal to $\epsilon (ba)$.  We obtain in particular the same type of
functorial
collection of flat bundles on smooth manifolds mapping to $X$, which gives a
flat bundle over any simplicial resolution of singularities as in \cite{Hodge3}.
The fundamental group of the topological realization of such a simplicial
resolution is the same as that of $X$ so this collection of flat bundles comes
from a local system on $X$. We call this the {\em underlying flat bundle} of the
variation ${\bf E}$.

\begin{theorem}
\label{TwistorIII}
{\rm (cf \cite{Hodge3} and  \cite{SteenbrinkZucker})}
If $X$ is a complex projective variety (possibly singular) with a Zariski open
subset $U\subset X$ and if
$(\Ee , W_{\cdot}, D)$ is a variation
of mixed twistor structure on $X$ with underlying flat bundle $V$, then the
cohomology $H^i(U, V|_U)$ carries a natural mixed twistor structure which is
functorial in $U$.
\end{theorem}

{\em Remark:} The same yoga of weights as in \cite{Hodge3}, Theor\`eme 8.2.4,
depending on the openness and singularity of the variety, holds here.

\subnumero{Mixed twistor complexes}
For the proofs of theorems \ref{TwistorII} and \ref{TwistorIII} (which we do
together) we proceed exactly as in \cite{Hodge2} and \cite{Hodge3}. In fact we
look at what is really going on in \cite{Hodge3} and plug in the stuff from
\cite{Hodge2}.

Start with the following definitions. A {\em mixed twistor complex}
is a filtered complex $(M^{\cdot}, W^{\rm pre}_{\cdot})$ of sheaves of $\Oo
_{\pp
^1}$-modules on $\pp ^1$ such that
$$
\underline{H}^i(Gr ^{W^{\rm pre}}_n(M^{\cdot}))
$$
(the cohomology sheaf) is a locally free sheaf of $\Oo _{\pp ^1}$-modules of
finite rank, pure of weight $n+i$. This shift of weights is the reason we
call this the {\em pre-weight filtration} and
use
the superfix $W^{\rm pre}$.

\begin{lemma}
\label{degen}
{\rm (\cite{Hodge3} Scholie 8.1.9)}
Suppose $(M^{\cdot}, W^{\rm pre}_{\cdot})$ is a mixed twistor complex. Then
the spectral sequence for a filtered complex which calculates the cohomology
sheaves  $\underline{H}^i(M^{\cdot})$ degenerates at $E_3$ (i.e. $d_r=0$ for
$r\geq 3$). The cohomology sheaves are locally free sheaves of $\Oo _{\pp
^1}$-modules and the filtration induced by the pre-weight filtration, when
shifted to $$
W_n\underline{H}^i(M^{\cdot}):= W^{\rm pre}_{n-i}\underline{H}^i(M^{\cdot})
$$
is the weight filtration for a mixed twistor structure on
$\underline{H}^i(M^{\cdot})$.
\end{lemma}
{\em Proof:}
The spectral sequence in question starts with
$$
E^{p,q}_1 = Gr ^{W^{\rm pre}}_{-p}(M^{p+q}),
$$
with differential $d_1: E_1^{p,q}\rightarrow E_1^{p,q+1}$ being the
associated-graded of the differential of $M^{\cdot}$.
The term $E^{p,q}_2$ is the cohomology of the differential $d_1$ and has
differential $d_2: E^{p,q}_2\rightarrow E^{p+1,q}_2$.
The cohomology of $d_1$
is just the cohomology sheaf which appears in the definition of mixed
twistor complex above (with $i=p+q$ and $n=-p$); thus, by definition $E^{p,q}_2$
is a locally free sheaf of $\Oo _{\pp ^1}$-modules, pure of weight $q$.  The
differential $d_2$ is a morphism of semistable bundles of the same slope
$q$, so
$E^{p,q}_3$, i.e. the cohomology of $d_2$, is again a locally free sheaf of $\Oo
_{\pp ^1}$-modules, pure of weight $q$. Now for $r\geq 3$ we have
$d_r: E^{p,q}_r\rightarrow E^{p+r-1, q+2-r}$ which is a morphism between
bundles which are, argueing inductively on $r$, pure of weights $q$ and $q+2-r <
q$ respectively; thus $d_r=0$ and the spectral sequence degenerates.

The limit of the spectral sequence is
$$
E^{p,q}_3 = Gr ^{W^{\rm pre}}_{-p}(\underline{H}^{p+q}(M^{\cdot})),
$$
and the fact that the associated graded pieces are locally free implies that
$\underline{H}^i(M^{\cdot})$ is locally free and its induced filtration is
by strict subbundles. The above $E^{p,q}_3$ being pure of weight $q = (p+q) +
(-p)$ we can rewrite as saying that
$W^{\rm pre}_{n-i}\underline{H}^i(M^{\cdot})$ is pure of weight $i+ (n-i)=n$.
This shows that the weight filtration as shifted in the statement of the lemma
gives a mixed twistor structure.
\eop

Our mixed twistor complexes will generally come from a more global situation.
Suppose $Y$ is a topological space mapping to $\pp ^1$ with map denoted $p$.   A
{\em mixed twistor complex on $Y/\pp ^1$} is a filtered complex of sheaves
$(N^{\cdot}, W^{\rm pre})$ of $p^{-1}(\Oo _{\pp ^1})$-modules on $Y$ such that
$$
{\bf R}^i p_{\ast} (Gr ^{W^{\rm pre}}_n(N^{\cdot}))
$$
is a locally free sheaf of $\Oo _{\pp ^1}$-modules, pure of weight $n+i$.
In this case, let ${\rm Go}(N^{\cdot})$ be the canonical Godement resolution
(or any other canonical flasque resolution), commuting with subquotients. Set
$$
M^{\cdot} = p_{\ast}{\rm Go}(N^{\cdot})
$$
with filtration $W^{\rm pre}M = p_{\ast}(W^{\rm pre}N^{\cdot})$.
The direct images of these resolutions calculate the higher direct images of the
complexes in question. Since $Gr ^{{\rm Go}(W^{\em pre})}_n({\rm Go}
(N^{\cdot})= {\rm Go} (Gr
^{W^{\rm pre}}_n(N^{\cdot})$  we have
$$
\underline{H}^i(Gr^{W^{\rm pre}M}_n(M^{\cdot})) =
{\bf R}^i p_{\ast} (Gr ^{W^{\rm pre}}_n(N^{\cdot}))
$$
which is, by hypothesis, pure of weight $n+i$. Thus $(M^{\cdot}, W^{\rm pre}M)$
is a mixed twistor complex on $\pp ^1$.

Of course in the case $Y=\pp ^1$ we recover the usual notion of mixed twistor
complex (noting however that in the passage from $N$ to $M$ in this case, the
complex will be replaced by its Godement resolution).

\subnumero{Patching}

The complexes constructed in \cite{Hodge2} for calculating the cohomology of
open varieties are closely tied to the holomorphic structure of $X$.  This
construction will only work over the standard ${\bf A}^1$---neighborhood of $0$.
A similar construction will work for the holomorphic structure $\overline{X}$
for the other standard neighborhood ${\bf A}^1$ of $\infty$. We need to patch
these two together, and here is where the relationship with the Leray spectral
sequence (pointed out by N. Katz, according to \cite{Hodge2}) comes in. Over
$\Gm$ this allows us to relate either of the previous two constructions with a
topological construction for flat bundles. We return to the trichotomy
$(\Ff , \nabla )$---$\{ L_t\}$---$(\Ff ', \nabla ')$ that was discussed in \S 3.

For each of these regimes we obtain a complex of sheaves on the appropriate open
subset of $\pp ^1$. We then need to  patch them together to obtain a mixed
twistor complex. This is the analogue in the present situation of the two
components, holomorphic and topological, in Deligne's definition of mixed Hodge
complex \cite{Hodge3} (our third component
$(\Ff ', \nabla ')$ being, in that case, the complex
conjugate of the first component
$(\Ff , \nabla )$, because Deligne considered real objects rather than
complex ones as we do here).

So we discuss in general how to patch together complexes of sheaves to obtain
mixed twistor complexes.

The basic situation we will treat is the following one: suppose we have
filtered complexes $M^{\cdot}$ and $N^{\cdot}$ of sheaves of $\Oo$-modules
(with the pre-weight filtrations denoted $W^{\rm pre}M$ etc.), respectively over
the standard neighborhoods ${\bf A}^1$ of $0$ and $\infty$ in $\pp ^1$.
Suppose we
have a filtered complex $P^{\cdot}$ of sheaves of $\Oo$-modules on $\Gm$
(the intersection of the two affine lines).  Finally suppose that we have
filtered quasiisomorphisms
$$
M^{\cdot}|_{\Gm} \leftarrow P^{\cdot} \rightarrow N^{\cdot} |_{\Gm} .
$$
We construct a filtered complex of sheaves of $\Oo$-modules on $\pp ^1$ denoted
$Patch (M\leftarrow P \rightarrow N)$ (or just $Patch$ for short) in the
following way.
\newline
(1)\, Let $i$ denote one of the three inclusions ${\bf A}^1\hookrightarrow \pp
^1$ (twice) or $\Gm \hookrightarrow \pp ^1$. Let $Ri_{\rm ex}$ denote a fixed
functorial choice of right derived functor (a real functor from complexes to
complexes, not  just a functor on the derived category) of some extension
functor
for the inclusion $i$. For example this could be $i_{\ast} \circ {\rm Go}$, the
composition of direct image with the Godement resolution. We require that there
be fixed a functorial quasiisomorphism
$$
\Ff ^{\cdot}\rightarrow i^{\ast} Ri_{\rm ex}(\Ff ^{\cdot})
$$
(which is the case  notably in our example).
\newline
(2)\, Let $M^{\cdot} _{\rm ex}:=
Ri_{\rm ex}M^{\cdot}$ be this extension functor applied to the filtered complex
$M^{\cdot}$
(resp. $P^{\cdot}_{\rm ex}:=Ri_{\rm ex}P^{\cdot}$, $N^{\cdot}_{\rm ex}
:=Ri_{\rm ex}N^{\cdot}$ with $i$ being the appropriate
inclusions).  They are again filtered complexes of sheaves of $\Oo
$-modules.
\newline
(3)\, In general if $f:A^{\cdot} \rightarrow B^{\cdot}$ is a map of filtered
complexes (or complexes of sheaves) let $Cone (A\rightarrow B)$ denote the
filtered  complex (or complex of sheaves) defined as
$$
Cone(A\rightarrow B)^k
:= A^{k+1} \oplus B^k
$$
with differential equal to $d_A + d_B + f$.
\newline
(4)\, Now note that we have a morphism
of filtered complexes of sheaves of $\Oo$-modules
$$
P^{\cdot}_{\rm ex} \rightarrow
M^{\cdot}_{\rm ex}\oplus N^{\cdot}_{\rm ex}
$$
(put in a minus sign in one of the factors).
We define
$$
Patch( M \leftarrow P \rightarrow N):= Cone
(P^{\cdot}_{\rm ex} \rightarrow M^{\cdot}_{\rm ex}\oplus
N^{\cdot}_{\rm ex}).
$$
We claim that
$$
\underline{H}^i(Gr^{W^{\rm pre}}_n Patch( M \leftarrow P \rightarrow
N))
$$
is the sheaf of $\Oo$-modules on $\pp ^1$ obtained by glueing together
$\underline{H}^i(Gr^{W^{\rm pre}}_nM^{\cdot})$ over the first neighborhood ${\bf
A}^1$, with $\underline{H}^i(Gr^{W^{\rm pre}}_nN^{\cdot})$ over the second
neighborhood ${\bf A}^1$, via the isomorphisms of cohomology sheaves induced by
the filtered quasiisomorphisms
$$
M^{\cdot}|_{\Gm} \leftarrow P^{\cdot} \rightarrow N^{\cdot} |_{\Gm} .
$$
Suppose $U$ is a connected open set contained in the
first affine neighborhood. Then $U\cap {\bf A}^1$ (intersection with the
other affine neighborhood) is equal to $U\cap \Gm$ and  $P|_{U\cap
\Gm}\rightarrow N|_{U\cap \Gm}$ is a filtered quasiisomorphism. Thus
$$
P^{\cdot}_{\rm ex}|_U\rightarrow N^{\cdot}_{\rm ex}|_U
$$
is a filtered quasiisomorphism. In general if $A \rightarrow B$ is a filtered
quasiisomorphism of complexes of sheaves and $A\rightarrow C$ is any morphism
then $B \rightarrow Cone (A \rightarrow B\oplus C)$ is a filtered
quasiisomorphism. In our case this says that the composition
$$
M^{\cdot} |_U \rightarrow M^{\cdot}_{\rm ex} |_U
\rightarrow Patch (M \leftarrow P \rightarrow N) |_U
$$
is a filtered quasiisomorphism
(the first arrow is the filtered quasiisomorphism coming from our assumption in
(2) above). Hence the cohomology sheaf of the associated graded of
$Patch$ restricted to $U$ is  naturally isomorphic to that  of $M$. Similarly if
$U'$ is a neighborhood contained in the second ${\bf A}^1$ then the cohomology
sheaf is naturally isomorphic to that of $N$. Finally, if $U$ is contained in
$\Gm$ then the composition of these two isomorphisms is equal to that which
comes from the diagram $M\leftarrow P\rightarrow N$ (the minus sign we put in
before garanties that we are not missing a minus sign here!).

Thus in order to construct a mixed twistor complex $Patch$ we just need to have
all of the above data with the property that the cohomology sheaves of the
associated-graded for the pre-weight filtration, when glued together, become
pure of the right weights.

We can make a similar construction in any more complicated situation of a chain
of quasiisomorphisms.  For example the actual situation we will need to
consider is when we have the sequence of filtered quasiisomorphisms
of filtered complexes over $\Gm$,
$$
M |_{\Gm} \leftarrow P \rightarrow Q \leftarrow R \rightarrow N|_{\Gm} .
$$
This can be replaced by the sequence
$$
M |_{\Gm} \leftarrow Cone '(P \oplus R \rightarrow Q) \rightarrow N|_{\Gm}
$$
where $Cone '$ is the cone but shifted in such a way that it is normalized for
the first variable (the cone we defined above was normalized for the second
variable). Now we can directly apply the previous discussion and define
$$
Patch (M,P,Q,R,N):= Patch (M \leftarrow Cone '(P \oplus R \rightarrow Q)
\rightarrow N).
$$
One
amusing point to notice is that in this case there are two stray minus signs
which cancel out. This is probably quite lucky if one wants to look at the real
situation, where a single minus sign might be very painful.

Finally, the whole thing goes through equally well in the relative situation. If
$p:Y\rightarrow \pp  ^1$ is a morphism of topological spaces and if we have
filtered complexes of sheaves $M$ and $N$ on the two $p^{-1}({\bf A}^1)$ and a
filtered complex of sheaves $P$ on $p^{-1}(\Gm )$ with filtered
quasiisomorphisms
$$
M|_{p^{-1}(\Gm )} \leftarrow P \rightarrow N|_{p^{-1}(\Gm )}
$$
then we obtain a filtered complex $Patch (M \leftarrow P\rightarrow N)$ of
sheaves on $Y\times \pp ^1$, which has the effect of patching together the
cohomology sheaves of the complexes $M$ and $N$ along $P$.  We can then take the
direct image down to $\pp ^1$ as was described previously.
On the other hand, we could also take the direct images of $M$, $N$ and $P$ down
to $\pp ^1$ and then patch them together. The answers in these two cases are not
quite the same but are quasiisomorphic. For simplicity we choose the route of
first taking the direct image then patching. Thus, to be concrete, we obtain a
filtered complex of sheaves
$$
Patch ( R^ip_{\ast}(M,P,N)):=
Patch \left( R^ip_{\ast}(M)\leftarrow R^ip_{\ast}(P) \rightarrow
R^ip_{\ast}(N)\right)
$$
on $\pp ^1$, patching together the higher direct images (the higher direct
images being calculated by canonical flasque resolutions, for example).

There is a
spectral sequence
$$
R^i p_{\ast}Gr ^{W^{\rm pre}}_n(M) \Rightarrow Gr ^{W^{\rm pre}}_nR^ip_{\ast}M
$$
and the same for $P$ and $N$, each of these being a spectral sequence of sheaves
over the appropriate open subsets of $\pp ^1$.  Note that filtered
quasiisomorphisms induce isomorphisms of spectral sequences.  Thus we may patch
together the spectral sequences to obtain a spectral sequence (of sheaves
on $\pp
^1$) which we denote
$$
Patch (R^i p_{\ast}Gr ^{W^{\rm pre}}_n(M,P,N)) \Rightarrow Gr ^{W^{\rm
pre}}_nPatch ( R^ip_{\ast}(M,P,N)).
$$
The terminology on the right was defined above and similarly the notation on the
left is defined as
$$
Patch (R^i p_{\ast}Gr ^{W^{\rm pre}}_n(M,P,N)):=
$$
$$
Patch \left( R^i p_{\ast}Gr ^{W^{\rm pre}}_n(M)\leftarrow
R^i p_{\ast}Gr ^{W^{\rm pre}}_n(P)\rightarrow
R^i p_{\ast}Gr ^{W^{\rm pre}}_n(N) \right) .
$$
If the beginning term of the spectral sequence is a
locally free sheaf of $\Oo _{\pp ^1}$-modules pure of weight $n+i$ then the
spectral sequence degenerates after the next term, and the answer is again pure
of weight $n+i$ (cf the argument of Lemma \ref{degen}).

Thus in order to obtain a mixed twistor complex by patching, it suffices that
the resulting patch of the higher direct images of the associated graded,
$$
Patch \, \left( R^i p_{\ast}Gr ^{W^{\rm pre}}_n(M),
R^i p_{\ast}Gr ^{W^{\rm pre}}_n(P),R^i p_{\ast}Gr ^{W^{\rm pre}}_n(N)
\right)
$$
be a locally free sheaf of $\Oo _{\pp ^1}$-modules, pure of weight $n+i$.

A similar statement works in the situation of patching $M,P,Q,R,N$ (which is the
case we will use).

\subnumero{Logarithmic complexes for mixed twistor structures}
Suppose  now that $Z$ is a smooth variety with a normal crossings divisor $D$,
and that $U= Z-D$. Suppose  $E$ is a variation of mixed twistor structure
on $Z$,
which we will now look at in terms of the three weighted objects
$\Ff$ (with $t$-connection $\nabla$), $L_t$ for $t\in \Gm$, and $\Gg$
which was denoted $\Ff '$ previously.  These lie respectively over ${\bf A}^1$,
$\Gm$ and the other ${\bf A}^1$.

We obtain the following complexes on subsets of the topological space $Y:=
Z^{\rm
top}\times \pp ^1$.  First, $M^{\cdot}$ is the complex over ${\bf A}^1$ of
holomorphic logarithmic differentials with coefficients in $\Ff$,
$$
M^i = \xi \Omega ^i_Z(\log D) \otimes _{\Oo_{X\times {\bf A} ^1} } \Ff
$$
with differential coming from the $t$-connection $\nabla$ and weight filtration
combining the weight filtration of \cite{Hodge2} with the filtration of $\Ff$.
Here as below, the symbol $\xi$ refers to the operation described \cite{NAHT},
\cite{SantaCruz}, with respect to the Hodge filtration of
$\Omega ^i_X(\log D)$  (which in this case is the ``stupid'' filtration).
Second, $N^{\cdot}$ is the same thing for $\Gg$ on $\overline{Z}$ transported to
$Y$ via the isomorphism $\overline{Z}^{\rm top} \cong Y$.  Finally, we define
several filtered complexes $P$, $Q$ and $R$ on the inverse image of $\Gm$.
Recall that over $t\in \Gm$ we have the local system $L_t$ which is the local
system of flat sections of $(\Ff , t^{-1} \nabla )$ and similarly for $\Gg$. Let
${\cal A}^{\cdot}_U$ denote the sheaf of $\Cc ^{\infty}$ differential forms on
$U$ and let $j: U\rightarrow Z$ denote the inclusion.  We have filtered
quasiisomorphisms (cf \cite{Hodge2} p. 33, the map called $\beta$ plus a
discussion analogous to that of the top of page 33)
$$
(\xi \Omega ^{\cdot}_Z(\log D) \otimes _{\Oo_{X\times \Gm }} \Ff , \tau )
\rightarrow
(j_{\ast}\xi \Omega ^{\cdot}_U\otimes _{\Oo_{X\times \Gm }} \Ff, \tau )
$$
$$
\rightarrow
(\xi {\cal A}^{\cdot}_U\otimes _{\Oo _{X\times \Gm }}\Ff , \tau )
= (\xi {\cal A}^{\cdot} _U \otimes _{p^{-1}\Oo _{\Gm }} \Ll , \tau ).
$$
At the end note that the family $\Ll$ is considered as a locally free sheaf of
$p^{-1}\Oo _{\Gm}$-modules on $Y\times _{\pp ^1} \Gm$.

Here the filtrations $\tau$ are the ``intelligent'' truncations of the complexes
of differentials, tensored with the weight filtrations of $\Ff$ or $\Ll$.

On the other hand, we have (again see \cite{Hodge2} p. 33, the map called
$\alpha$) the filtered quasiisomorphism over $\Gm$,
$$
(\xi \Omega ^{\cdot}_Z(\log D) \otimes _{\Oo_{X\times \Gm}} \Ff , W^{\rm pre} )
\leftarrow
(\xi \Omega ^{\cdot}_Z(\log D) \otimes _{\Oo_{X\times \Gm}} \Ff , \tau ) .
$$
In view of this we put
$$
P^{\cdot} := (\xi \Omega ^{\cdot}_Z(\log D) \otimes _{\Oo_{X\times \Gm}} \Ff ,
\tau )
$$
over $\Gm$ and
$$
Q^{\cdot} := (\xi {\cal A}^{\cdot} _U \otimes _{p^{-1}\Oo _{\Gm}} \Ll , \tau )
$$
again over $\Gm$. Note that $Q$ now makes no reference to the complex structure
so it is the same with respect to $Z$ or $\overline{Z}$.  Let $R$ be the same
thing as $P$ but on $\overline{Z}$ and using $\Gg$.  We obtain the diagram of
filtered quasiisomorphisms of complexes on $Y\times _{\pp ^1} \Gm$,
$$
M |_{\Gm} \leftarrow P \rightarrow Q \leftarrow R \rightarrow N|_{\Gm} .
$$
Put
$$
MTC(E):= Patch (Rp_{\ast}M,Rp_{\ast}P,Rp_{\ast}Q,Rp_{\ast}R,Rp_{\ast}N)
$$
with these and the previous notations.
We claim that this is a mixed twistor complex. This claim implies that the
cohomology sheaves of $MTC(E)$ have natural mixed twistor structures. The fiber
over $1$ is the same as the hypercohomology on $Z^{\rm top}$ of the complex $Q$,
which is just the cohomology of $U$ with coefficients in the local system $L_1
|_U$.  This will therefore prove Theorem \ref{TwistorII}.

\subnumero{Proof of claim}
An {\em extension} of filtered complexes of sheaves is a short exact sequence
of complexes inducing short exact sequences on all levels of the filtrations.
This induces an extension of the associated-graded complexes.

An extension of VMTS
$$
0 \rightarrow E ' \rightarrow E \rightarrow E'' \rightarrow 0
$$
induces an extension of filtered complexes
$$
0 \rightarrow MTC(E ') \rightarrow MTC(E) \rightarrow MTC(E'') \rightarrow 0,
$$
hence an extension
$$
0 \rightarrow Gr _n^{W^{\rm pre}}MTC(E ') \rightarrow Gr _n^{W^{\rm pre}}
MTC(E) \rightarrow Gr _n^{W^{\rm pre}}MTC(E'')
\rightarrow 0.
$$
From this we get a long exact sequence of cohomology sheaves.
Suppose we know that
$$
\underline{H}^iGr _n^{W^{\rm pre}}MTC(E ') \;\; \mbox{and} \;\;
\underline{H}^iGr _n^{W^{\rm pre}}MTC(E '')
$$
are pure of weight $n+i$, and suppose we know also that the connecting maps
in the
long exact sequence are
zero.  Then we can conclude that $\underline{H}^iGr _n^{W^{\rm pre}}MTC(E )$
are pure of weight $n+i$. This last phrase may be restated as saying that if
we know the connecting maps in the long exact sequence of cohomology of
associated graded pieces are zero, then $MTC(E ')$ and $MTC(E '')$
being mixed twistor complexes implies that $MTC(E )$ is a mixed twistor complex.

We prove that for any VMTS $E$, $MTC(E)$ is a mixed twistor complex, by
induction on the size of the interval containing the weight-graded pieces of
$E$. If the size of this interval is $1$ then $E$ is pure and we will treat
this case below.  Suppose that $E$ is in an interval of size $n$ and that
we know
the result for any VMTS in an interval of size $<n$. Let $W_iE$ be the lowest
nonzero piece of the weight filtration. Apply the preceeding discussion to the
extension of VMTS  $$
0\rightarrow W_iE \rightarrow E \rightarrow E/W_iE \rightarrow 0.
$$
By our inductive hypothesis $MTC(W_iE)$ and $MTC(E/W_iE)$ are mixed twistor
complexes. In order to be able to conclude that $MTC(E)$ is a mixed twistor
complex, we just have to know that the connecting maps in the long exact
sequence calculating $H^jGr ^{W^{\rm pre}}_nMTC(E)$, are zero. But these
connecting maps are obtained by a construction of the form, take an element with
coefficients in $E'' := E/W_iE$ and lift it to an element with coefficients in
$E$, then take the coboundary which will have coefficients in $E' := W_iE$.
The preweight of our element (which lies in some complex of forms or Godement
resolution) is the sum of the weight of the coefficient in
$E$ plus the weights of the other stuff.  But the weight of the coefficient in
$E$ is strictly decreased by the operation described above, since the weights of
$E'$ are strictly lower than the weights of $E''$. Thus, this operation is zero
on the associated graded for the preweight filtration, in other words the
connecting map in question is zero.

This completes the proof of the claim modulo the case where $E$ is pure, which
we now treat. It suffices to consider the case where $E$ is pure of weight
zero.
By the discussion of the previous subsection (and with the same notations as
there), it suffices to check that the patching of $R^ip_{\ast}Gr^{W^{\rm
pre}}_n({\bf x} )$ for ${\bf x} = M,P,Q,R,N$, be pure of weight $n+i$.

Let $D^{(k)}$ denote the disjoint union of the intersections of $k$ smooth
components of $D$.  We assume that the irreducible components of $D$ are smooth
so that the sheaves denoted $\varepsilon$ in \cite{Hodge2} are trivial.

The {\em residue map} gives an isomorphism
$$
res: Gr^{W^{\rm pre}}_n(M^{i}) \cong \xi \Omega ^{i-n}_{D^{(n)}}\otimes \Ff
|_{D^{(n)}} \otimes \Oo _{{\bf A}^1}(n\cdot 0)
$$
with differential induced by the $t$-connection $\nabla$.
The term $\Oo _{{\bf A}^1}(n\cdot 0)$ (which means the sheaf of functions with
poles of order $n$ at the origin) comes from the fact that the residue map
contracts out $n$ things of the form $dz_i /z_i$, which provide twists due to
the construction $\xi$.
Over
$t\neq 0$ (i.e. over $\Gm \subset \pp ^1$) this complex is a resolution of the
local system $L_t |_{D^{(n)}}$.  Thus we have
$$
R^ip_{\ast}Gr^{W^{\rm pre}}_n(M^{\cdot}) =
R^{i-n}p_{\ast}(\xi \Omega ^{\cdot}_{D^{(n)}}\otimes \Ff
|_{D^{(n)}}) \otimes \Oo _{{\bf A}^1}(n\cdot 0),
$$
and this restricts (via the natural isomorphism) to $R^{i-n}p_{\ast}(\Ll
|_{D^{(n)}})$ over $\Gm$.  The direct images of the associated-graded pieces of
the complexes $P,Q,R$ all give the same answer $R^{i-n}p_{\ast}(\Ll
|_{D^{(n)}})$.

The same holds for $R^{i-n}p_{\ast}Gr^{W^{\rm pre}}_n(N^{\cdot})$ and again
this restricts to $R^{i-n}p_{\ast}(\Ll |_{D^{(n)}})$ (which doesn't depend
on the
complex structure of $D^{(n)}$).

The patching of these direct images of associated-graded pieces for $M,P,Q,R,N$
yields the filtered coherent sheaf constructed in Theorem \ref{coho1}, for our
VMTS $E$ pulled back to $D^{(n)}$, tensored by
$\Oo _{\pp ^1}(2n)$ because of term $\Oo _{{\bf A}^1}(n\cdot 0)$ and the similar
term in the  neighborhood of $\infty$.    By the result of \S 4 (recall
that here
we are treating the case where $E$ is pure of weight $0$) this patching is pure
of weight $(i-n)+ 2n =n+i$. This completes the proof of the claim (and hence the
proof of Theorem \ref{TwistorII}).

\subnumero{The simplicial situation}

In order to prove Theorem \ref{TwistorIII} we need to consider a simplicial
situation. Any simplicial scheme $X_{\cdot}$ (of finite type) can  be replaced
(via the method of \cite{Hodge3}) by a simplicial collection of projective
smooth
schemes $Z_{\cdot}$ containing normal crossing divisors $D_{\cdot}$ such
that the
maps are compatible with the normal crossing divisors. The topological type of
the  original simplicial scheme is recovered by taking the topological type
of the
simplicial open subschemes $U_{\cdot} = X_{\cdot} - D_{\cdot}$ (complement
of the
divisor with normal crossings).

For each $k$ we have constructed a mixed twistor complex $(M^{\cdot}_k, W^{\rm
pre})$ calculating the cohomology of $U_k$ using the normal crossings
compactification $(X_k, D_k)$.  This is functorial (in a contravariant way): the
face maps give morphisms of mixed twistor complexes
$$
(M^{\cdot}_k, W^{\rm
pre}) \rightarrow (M^{\cdot}_{k+1}, W^{\rm
pre}).
$$
In other words, we get a {\em cosimplicial mixed twistor complex}.  The
cohomology of the simplicial scheme is calculated by taking the associated
double complex (with new differential the alternating sum of the face maps)
and then taking its associated single complex.  We just have to
provide this with a structure of mixed twistor complex. This is of course
exactly what is explained in \cite{Hodge3}. We briefly sketch the argument.
We obtain a total
complex $$
N^j = \bigoplus _{i+k=j}M^i_k,
$$
and we define its preweight filtration by
$$
W^{\rm pre}_nN^j := \bigoplus _{i+k=j} W^{\rm pre}_{n+k}M^i_k.
$$
The two differentials are compatible with this new weight filtration
and in fact the differential $M^i_k \rightarrow M^i_{k+1}$ induces
the zero map on the associated-graded for the weight filtration. In particular,
the double complex $Gr^{W^{\rm pre}}_n(M^{\cdot}_{\cdot})$ has one of its
differentials vanishing; the other differential is just the differential in each
mixed twistor complex $M^{\cdot}_k$.  Thus
$$
H^j(Gr^{W^{\rm pre}}_n(M^{\cdot}_{\cdot}))= \bigoplus _k H^{j-k}(Gr^{W^{\rm
pre}}_{n+k}(M^{\cdot}_k)),
$$
and this is pure of weight $n+j$ since each of the $M^{\cdot}_k$ is a mixed
twistor complex. This shows that the total complex is again a mixed twistor
complex.

The proof that the resulting mixed twistor structure on the cohomology is
independant of the choice of desingularization and normal crossing
compactification, (and, what is pretty much the same thing, functoriality) is
exactly the same as in \cite{Hodge2} \cite{Hodge3}.
This completes the proof of Theorem \ref{TwistorIII}.
\eop

We have also obtained the
simplicial version:

\begin{theorem}
\label{TwistorSimpl}
Suppose $X_{\cdot}$ is a simplicial scheme with each $X_k$ projective over $Spec
(\cc )$, and suppose $U_{\cdot} \subset X_{\cdot}$ is a simplicial open
subscheme. Suppose $(E, W, {\bf d})$ is a graded-polarizable variation of mixed
twistor structure on $X_{\cdot}$ (i.e. a functorial association of VMTS for
every
smooth scheme mapping to $X_k$, with isomorphisms of compatibility for the
simplicial maps). Then there is a mixed twistor structure on $H^i(U_{\cdot}, L)$
where $L$ is the flat bundle associated to $(E_1, {\bf d}_1)$.
This mixed twistor structure is functorial for maps of VMTS and also for maps of
simplicial schemes $U_{\cdot}$ (from the cohomology of $E$ on the target to the
cohomology of the pullback of $E$ on the domain).
\end{theorem}
\eop

{\em Remark:}  For now, we have not made any study of VMTS on open varieties,
so we are forced to assume the existence of extensions to projective varieties
$X_k$ and furthermore that these extensions should be organized in a simplicial
way. Eventually these restrictions should be gotten rid of (possibly in
exchange for a condition of admissibility of the VMTS on the open varieties).

\numero{Nilpotent orbits and the limiting mixed twistor structure}

Suppose $U=X-P$ is a compact curve with one point removed (actually any number
of points is ok but we work near one) and suppose $E$ is a harmonic bundle on
$U$. We suppose that the monodromy of $E$ at $P$ is unipotent, and that the
flat sections have sub-polynomial growth. This last condition means that, in the
language of \cite{HBNC}, the associated filtered local system has trivial
filtration.

Let $\Ee$ be the associated variation of twistor structure considered as a
bundle over $U\times \pp ^1$, normalized to have weight zero say.  Let $(\Ff _U,
\Ll , \Gg _U)$ be the triple of a $\lambda$-connection $\Ff_U$, a family of
local
systems $\Ll$ on $U$, and a $\lambda ^{-1}$-connection $\Gg _U$ on
$\overline{Z}$. Let $\Ff$ be the canonical
extension of $\Ff_U$ to a logarithmic $\lambda$-connection on $Z$. Let $\Gg$ be
the canonical extension of $\Gg_U$ to a logarithmic $\lambda
^{-1}$-connection on
$\overline{Z}$. These extensions have nilpotent transformations in the fiber
over $P$, hence they have monodromy weight filtrations (cf \cite{Schmid} for
example).  The dimensions of the associate-graded of the weight filtrations are
the same for all $\lambda $ (cf \cite{HBNC}). This implies that the weight
filtrations are filtrations of the bundles $\Ff _P$ (defined to be the
restriction to ${\bf A}^1 = \{ P\} \times {\bf A}^1$) and $\Gg _P$ (similarly,
the restriction to the other copy of $\{ P\} \times {\bf A}^1$), by strict
subbundles.

Fix a tangent vector $\eta$ at $P$, fix  a point $P_0$ nearby $P$, and fix a
path going from $P_0$ to $P$ arriving by the direction $\eta$.  For any
logarithmic flat connection on a neighborhood of $P$ there is an isomorphism
with a fixed connection with constant coefficients on a trivial bundle, and this
isomorphism determines an isomorphism between the fiber over $P$ and the fiber
over $P_0$.  The isomorphism of fibers is well defined depending on the choice
of tangent vector \cite{DeligneRegSing}.  Furthermore this
isomorphism depends analytically on any parameters, and takes the monodromy
operator on the fiber over $P_0$ to the exponential of the residue of the
connection over $P$.  In particular it preserves the weight filtrations.

Applying this to our situation, note that the family $\Ll _{P_0}$ of fibers of
the $L_t$ over $P_0$ has a family of weight filtrations. The above isomorphisms
give an isomorphism of filtered bundles
$$
\Ll _{P_0} \cong \Ff _P
$$
over $\Gm$.  Similarly we have an isomorphism of filtered bundles
$$
\Ll _{P_0} \cong \Gg _P
$$
over $\Gm$ and we can use these to glue together the filtered bundles $\Ff_P$
and  $\Gg_P$ to give a filtered bundle $(M, W)$ on $\pp ^1$.

Recall that $T(1)$ denotes the Tate twistor, of rank one and weight two.
It is isomorphic to the tangent bundle of $\pp ^1$ (and this is probably the
right way to say canonically what it is).

\begin{conjecture}
\label{LimitingMTS}
{\rm (\cite{Schmid})}
The filtered bundle defined above is a mixed twistor structure. (We call it the
{\em limiting mixed twistor structure of $E$}).  Furthermore the logarithms of
monodromy operators glue together to give a morphism of mixed twistor structures
$$
N: M \rightarrow M\otimes _{\Oo _{\pp ^1}}T(1)
$$
such that the weight filtration is deduced from $N$ in the now standard way.
\end{conjecture}

Once this is proved, we would then like to generalize the nilpotent orbit
theorem \cite{Schmid}, saying that the asymptotic behavior of our variation of
twistor structure is well approximated by the behavior of a standard variation
(``nilpotent orbit'') deduced directly from the limiting mixed twistor
structure.

All of this should then be generalized to higher dimensions, where it might help
us to get a hold on the behavior of harmonic bundles near normal crossings
divisors.

\numero{Jet-bundles of hypercomplex manifolds}

Let $U$ be a smooth hyperk\"ahler or more generally hypercomplex  manifold.
\footnote{
The definition of {\em hypercomplex manifold} is essentially the
that of {\em hyperk\"ahler manifold} (cf \cite{HKLR})  minus the condition of
existence of a metric.
The original reference for hypercomplex manifolds seems to be
Boyer \cite{Boyer}---see also the preprint of Kaledin \cite{Kaledin}.
}
The product
$T:= U \times \pp ^1$ has a complex structure called the {\em twistor space} of
the quaternionic structure of $U$ \cite{Hitchin1} \cite{Hitchin2} \cite{HKLR}
\cite{Kaledin}. Basically an almost-hypercomplex structure is just a $\Cc
^{\infty}$ family of actions of the quaternions on the tangent spaces of $U$,
which gives an almost-complex structure on $T$. Integrability of this complex
structure is exactly the condition of $U$ being hypercomplex (see Theorem 1 of
\cite{Kaledin}, or \cite{HKLR} for the hyperk\"ahler case). The projection
$p:T\rightarrow \pp ^1$ and the horizontal sections for the product structure
(called {\em twistor lines} \cite{HKLR}) $\zeta : \pp ^1\rightarrow T$ are
holomorphic maps. Fix a twistor line $\zeta$ corresponding to a point $u\in U$.
The normal bundle $N_{\zeta /T}$---a holomorphic vector bundle on $\pp ^1$ of
rank $r$ equal to $\frac{1}{2}dim _{\rr}(U)$---is the twistor-bundle of the
tangent space $T(U)_u$ with its quaternionic structure. In particular, $N_{\zeta
/T}$ is a pure twistor structure of weight $1$, i.e. it is isomorphic to a
direct
sum of  $r$ copies of $\Oo _{\pp ^1}(1)$.

Let $J^n_{\zeta /T}$ denote the $n$-th relative jet bundle of $T$ over
$\pp ^1$ along the section $\zeta$.  It may be described as follows. Let
$I_{\zeta/ T}$ denote the ideal sheaf of the section $\zeta$. The coherent
sheaf $Q_{n, \zeta /T} := \Oo _T / I_{\zeta / T} ^{n+1}$ on $T$ is flat and
proper
over $\pp ^1$. Its direct image $p _{\ast}(Q_{n, \zeta /T} )$ is a vector
bundle on $\pp ^1$
and the jet bundle is the dual of this vector bundle.  We have inclusions
$$
J^0 _{\zeta /T} \subset \ldots \subset J^n_{\zeta /T}
$$
and the quotients are naturally symmetric powers of the normal bundle:
$$
J^m_{\zeta /T} / J^{m-1}_{\zeta /T} \cong Sym ^m(N_{\zeta /T}).
$$
In particular, these quotients are pure twistor structures of weight $m$.
Thus if we put $W_m(J^n_{\zeta /T}) := J^m_{\zeta /T}$ we obtain a mixed
twistor structure. The fiber over $1$ is the jet space $J^n_u$ of $U$
at $u$ (for the complex structure corresponding to $1\in \pp ^1$---this
depends on
how we normalize the twistor space). We conclude the phrase ``the jet spaces of
a hypercomplex manifold have natural mixed twistor structures''.

Taking the direct limit or union of the $J^n_u$ we get an ind-object
$J^{\infty}_u$ in the category of mixed twistor structures. The weight quotients
are finite dimensional, and in fact the finite $J^n_u$ are recovered as the
pieces of the weight filtration of $J^{\infty}_u$. By abuse of notation we will
say that $J^{\infty}_u$ is a mixed twistor structure which determines in
particular the mixed twistor structures $J^n_u$.

The dual of the jet space $J^n_u$ is $Q_{n,u}:=\Oo _{U,u} /{\bf
m}_u^{n+1}$. This
also has a mixed twistor structure, dual of the above.  The algebra structure of
$Q_{n,u}$ comes from a coalgebra structure of
$J^{\infty}_u$ (coalgebra in the category of ind-mixed twistor structures).
The mixed twistor structure $J^{\infty}_u$ together with its coalgebra
structure determines the algebras $Q_{n,u}$. In fact the
twistor structure gives the coherent sheaves of algebras denoted $\Ff$ above,
and taking the $Spec$ at each stage and formally taking the direct limit we
obtain the formal completion of $T$ along the section $\zeta$.  Thus
``the mixed twistor structure on the jet space $J^{\infty}_u$ of $U$ at $u$
determines the formal completion of the twistor space $T$''.

Finally the twistor space $T$ has an antilinear involution $\sigma _T$
covering the antipodal involution $\sigma _{\pp ^1}$. In particular the mixed
twistor structure $J^{\infty}_u$ has an antipodal real structure and, via the
previous paragraph, from this real structure we recover the involution $\sigma
_T$ on the formal completion along the section $\zeta$. Nearby the section
$\zeta$ the only other sections invariant under the involution $\sigma$ are the
twistor lines (this local uniqueness of solutions of an elliptic equation
may be verified in the linearization which is the normal bundle $N_{\zeta /T}$).
Hence the antipodal real structure of $J^{\infty}_u$ determines the product
structure of the formal completion of $T$. Therefore it determines the
hypercomplex structure of $U$ in the formal neighborhood of $u$.

{\em Question:} How can we include the data of the hyperk\"ahler metric in this
point of view?

We obtain a dictionary between formal germs of integrable
quaternionic structure $(u, \hat{U})$ and mixed twistor structures
$J^{\infty}$ with coalgebra structure and antipodal real structure (these last
two being compatible) such that the associated graded of the weight
filtration is
the symmetric algebra. This gives a coordinate-free description of the power
series of a quaternionic structure.

\subnumero{The ``Gauss map''}
Suppose $U$ is a hypercomplex manifold of real dimension $4d$, and fix a
positive integer $n$. Then for each $u\in U$ we have constructed a mixed twistor
structure on $J^n_u(U)$. The dimensions of the weight-graded quotients are
$(a_0, a_1, \ldots , a_n)$ with $a_i$ being the rank of $Sym ^i (\cc ^{2d})$.
We obtain a morphism to the moduli stack of mixed twistor structures
$$
U\rightarrow {\cal MTS}(a_0,\ldots , a_n).
$$
We can think of this as a ``Gauss map'' for the hypercomplex structure.
Actually
this can be lifted a little bit more, using the coalgebra structure on the jet
bundles.  We have an isomorphism of complex vector spaces (using the complex
structure $j$, say)
$$
Gr ^W_m(J^n_u)\cong Sym ^{m}(T(U)_u).
$$
In particular, if we choose a basis for $T(U)_u$ then this gives a basis for
each of the $Gr ^W_m(J^n_u)$, in other words a framing. Thus, over the frame
bundle of the tangent bundle we obtain a lifting of our map to
$Fr{\cal MTS}(a_0,\ldots , a_n)$.

We can be a little bit more precise about
this. Note that the associated-graded of the tensor product of two filtered
vector spaces, is naturally isomorphic to the tensor product of the
associated-graded vector spaces. In this way, a coalgebra structure on a
filtered vector space $J$ (which we now assume infinite dimensional, no longer
fixing $n$) gives a coalgebra structure on $Gr^W(J)$. If $T$ is a vector space
then let $Fr{\cal MTS}Cog(Sym^{\cdot}T)$ denote the moduli stack of
pairs $(J, \beta )$ where $J$ is a mixed twistor structure with coalgebra
structure and $\beta : Gr ^W(J_1)\cong Sym^{\cdot}T$ is an isomorphism of
coalgebras. But this  moduli stack is in fact a projective limit of schemes,
because of the rigidity property of the remark following Lemma \ref{abelian}.
Note also that  the coalgebra structure, if it exists, is uniquely determined by
$\beta$ again by the same rigidity property.  Thus letting $a_i$ denote the
dimension of $Sym ^i(T)$ we have
$$
Fr{\cal MTS}Cog(Sym^{\cdot}T) \subset \lim _{\leftarrow , n}
Fr{\cal MTS} (a_0,\ldots , a_n).
$$
Finally, $GL(T)\cong GL(a_1,\cc )$
acts on $Fr{\cal MTS}Cog(Sym^{\cdot}T)$ and our {\em Gauss map} is a
real analytic map
$$
\Phi : U\rightarrow Fr{\cal MTS}Cog(Sym^{\cdot}T)/GL(T).
$$
The associated $GL(T)$-bundle on $U$ is the tangent bundle (for the complex
structure over $1\in \pp ^1$).

{\em Problem:} what are the differential equations satisfied by $\Phi$
corresponding to the fact that it parametrizes the family of jet spaces of a
manifold $U$?

\numero{The moduli space of representations}
Suppose $X$ is a connected smooth projective variety with basepoint $x$. Let
$G=GL(n)$ (or any other reductive group) and let $M(X,G)$ denote the moduli
space
of representations of $\pi _1(X, x)$ in $G$ up to conjugacy
\cite{Lubotsky-Magid}.
The set of smooth points of $M(X,G)$ has a hyperk\"ahler structure
\cite{Hitchin1}
\cite{Hitchin2} \cite{SantaCruz}.  Thus we can apply the above remarks to any
smooth open set $U\subset M (X, G)$. As remarked in (\cite{SantaCruz} \S 3), any
smooth stratum for a canonical stratification of $M(X,G)$ also inherits a
hyperk\"ahler structure so we could also apply the previous remarks there.  We
find that the jet spaces $J^n_{\rho}(M(X,G))$ at any smooth point $\rho$ have
natural mixed twistor structures, and the same for jet spaces of smooth strata.

This can be generalized to singular points and points of the representation
spaces as follows.  In general if $Z$ is a scheme and $z\in Z$, denote by
$Q_{n,z}(Z)$ the algebra $\Oo _{Z,z}/{\bf m}_z^{n+1}$, and let $J^n_z(Z)$ denote
the $\cc$-linear dual coalgebra.  For our pointed connected smooth projective
variety $(X,x)$ as above, let $R(X,G,x)$ denote the scheme of representations of
$\pi _1(X,x)$ in $G$ (before dividing out by conjugation). Note that
$M(X,G)$ is the geometric invariant theory categorical quotient of $R(X,G,x)$ by
the action of $G$. Deligne's construction as described in \cite{NAHT}
\cite{SantaCruz}  gives an analogue of the twistor space for $R(X,G,x)$,
namely a
complex analytic space $R_{Del}(X,G,x)$ over $\pp ^ 1$ with ``twistor lines''
corresponding to the semisimple representations of $\pi _1(X,x)$ in $G$. The
group $G$ acts and the categorical quotient denoted $M_{Del}(X,G)$ gives, over
the smooth points, the twistor space for $M(X,G)$.  In Deligne's interpretation,
$M_{Del}(X,G)$ constitutes a twistor space for the ``singular hyperk\"ahler
structure" on $M(X,G)$ (in particular there is an antipodal involution
$\sigma$).

Suppose $\eta : \pp ^1 \rightarrow R_{Del}(X,G,x)$ is a twistor line
corresponding to a semisimple representation $\rho$. By Goldman-Millson theory
\cite{Goldman-Millson} (cf \cite{SantaCruz} and \cite{Moduli} for the present
application) the scheme $R_{Del}(X,G,x)$ is analytically formally trivial along
the section $\eta$ and furthermore the singularities in the transversal
direction are quadratic.  In particular we can form the family of jet bundles
$J^n_{\eta} (R_{Del}(X,G,x))$ which is a bundle of coalgebras over $\pp ^1$ (and
which we denote simply by $J^n$ for short, below). This bundle is dual to the
bundle of algebras $Q_{n, \zeta}(R_{Del}(X,G,x))$.

\begin{theorem}
\label{JetBdlRepSpace}
The jet bundles $J^n_{\eta} (R_{Del}(X,G,x))$ have natural mixed twistor
structures.
\end{theorem}
{\em Proof:}
The group $G$ acts.  Let $G\eta$ denote the orbit (fiber-by-fiber) of the
section $\eta$.
Let $W_kQ_{n, \zeta}(R_{Del}(X,G,x))$ be the $k$-th power of the ideal of
$G\eta \cap Spec (
Q_{n, \zeta}(R_{Del}(X,G,x)))$ in $Q_{n, \zeta}(R_{Del}(X,G,x))$.
This is the weight filtration.

By Luna's etale slice theorem
the space $R_{Del}(X,G,x)$ has a product structure formally around the section
$\eta$,
$$
R_{Del}(X,G,x)^{\wedge} = G\eta ^{\wedge} \times Z^{\wedge} .
$$
(to be precise this holds locally in the etale topology over the base $\pp ^1$).

We will use a different filtration to study this weight filtration.  Define the
{\em tangent cone filtration} $V_k$ of $Q_{n,u} (U)$ to be the filtration by the
powers of the maximal ideal. To fix notations, let
$$
Q_{\infty ,u}(U) = \lim _{\leftarrow , n} Q_{n ,u}(U)
$$
and put
$$
V_{-k}Q_{\infty ,u}(U):= {\bf m}_u ^kQ_{\infty ,u}(U).
$$
 The terminology comes from the
fact that $Spec (Gr ^V(Q_{\infty,u}(U))$ is the  tangent cone of $U$ at $u$.
Note of course that
$$
Q_{n ,u}(U) = Q_{\infty ,u}(U)/V_{-n-1}Q_{\infty ,u}(U).
$$
The dual of $V$ is a filtration which we also denote by $V_k$ (indexed in the
positive range of values of $k$ this time) of the jet space $J^{\infty}_u(U)$.

The above discussion works in our family (since the families are locally
trivial over $\pp ^1$) to give filtrations $V$ of $Q_{\infty , \eta}
(R_{Del}(X,G,x))$ and $J^{\infty}_{\eta}(R_{Del}(X,G,x))$.

If $U$ is a product, $U= U_1 \times U_2$ then we have
$$
Gr ^V(Q_{\infty , u}(U))= Gr ^V(Q_{\infty , u}(U_1))\otimes Gr ^V(Q_{\infty ,
u}(U_2)),
$$
in other words, the tangent cone is the product of the two tangent
cones. This is proved by interpreting the tangent cone as the scheme cut out by
the initial forms of the equations.

Now suppose that $U$ is a product as above and
let $W$ be the filtration by powers of the ideal of the factor $U_1$.
The filtration induced by $W$ on $Gr ^V(Q_{\infty , u}(U))$ is again the
filtration by powers of the ideal of the first factor. This is the same as the
filtration induced by the grading of the second factor only, in the above tensor
product structure.

By the lemma of two
filtrations (\cite{Hodge2}) we have
$$
Gr ^W_m Gr ^V_k (Q_{\infty , u}(U))=
Gr ^V_k Gr ^W_m (Q_{\infty , u}(U)).
$$
The same holds  in the situation of our family over $\pp^1$. Thus in order to
prove that $Q_{\infty , \eta}(R_{Del}(X,G,x))$ is a mixed twistor structure, it
suffices to prove that
$$
Gr ^V_{\cdot}Q_{\infty , \eta}(R_{Del}(X,G,x))
$$
is a mixed twistor structure.

Note that
$$
Gr ^V_{k}Q_{\infty , \eta}(R_{Del}(X,G,x))=
(J^k/J^{k-1})^{\ast},
$$
and in particular
$$
Gr ^V_{1}Q_{\infty , \eta}(R_{Del}(X,G,x))=
(J^1/J^{0})^{\ast} = N^{\ast}
$$
where $N$ denotes the normal bundle (isomorphic to $J^1/J^0$).  The algebra
structure of $Gr^V_{\cdot}$ gives a morphism
$$
Sym ^{k}(N^{\ast}) \rightarrow
(J^k/J^{k-1})^{\ast} \rightarrow 0
$$
which is a surjection because $V$ is by definition the filtration by powers
of the maximal ideal.  By Goldman and Millson-theory the family of algebras
$Q_{\infty , \eta}(R_{Del}(X,G,x))$ is quadratic, hence isomorphic to the
tangent
cone and in particular the tangent cone is quadratic. This means that the kernel
of the above surjection comes from the quadratic relations, in other words we
have an exact sequence
$$
P\otimes  Sym ^{k-2} (N^{\ast}) \rightarrow Sym ^{k}(N^{\ast}) \rightarrow
(J^k/J^{k-1})^{\ast} \rightarrow 0
$$
where $P \subset
Sym ^2(N^{\ast})$ is the module of quadratic relations, that is the kernel of
the map
$Sym ^2(N^{\ast}) \rightarrow J^2 / J^1$.

We claim that the filtration induced by $W$ on $(J^k/J^{k-1})^{\ast}$
is the image of the filtration of $Sym ^{k}(N^{\ast})$ induced by $W$ on
$N^{\ast}= (J^1/J^0)^{\ast}$. This is a general fact coming from the product
situation $U=U_1\times U_2$ with our filtrations $V$ and $W$ as above.
Write $U_i= Spec (A_i)$ and let $\germ_i\subset A_i$ be the maximal ideal of
the origin. Then $U_1\times U_2 = Spec (A)$ with $A=A_1\otimes A_2$,
$W$ is the filtration
by powers of $\germ _2A$ and $V$ is the filtration by powers of $\germ_1 A+
\germ_2A$.  The statement which needs to be proved $(\ast )$ is that
$$
(\germ ^1A + \germ ^2A) ^{\otimes p} \otimes (\germ _2A)^{\otimes q}
\rightarrow
\frac{(\germ _2A)^q \cap (\germ _1 A+ \germ _2A)^{p+q}}{(\germ _1 A+
\germ _2A)^{p+q+1}}
$$
is surjective. Choose a splittings $A_i\cong Gr (A_i)$ which we think of as
direct sum decompositions $A_i \cong \bigoplus A_i^k$ compatible with the
filtrations by powers of the maximal ideal. The product structure has an upper
diagonal form for this decomposition, with the product in $Gr(A_i)$ along the
diagonal.  We get a decomposition $A= \bigoplus A^{j,k}$ with $A^{j,k}=
A_1^j\otimes A_2^k$. The decomposition of $A$ into pieces of the form
$\bigoplus _{j+k=n}A^{j,k}$ corresponds to a splitting $A\cong Gr(A)$.
In particular,
$$
(\germ _1 A+ \germ _2A)^{n} = \bigoplus _{j+k\geq n}A^{j,k}.
$$
Similarly
$$
(\germ _2A)^q= \bigoplus _{k\geq q}A^{j,k}.
$$
From this we get
$$
(\germ _2A)^q \cap (\germ _1 A+ \germ _2A)^{p+q} = \bigoplus _{j+k \geq n, k\geq
q}A^{j,k} .
$$
The product on the associated graded induces a surjective morphism
$$
(A^{1,0})^{\otimes j}\otimes (A^{0,1})^{\otimes k}
\rightarrow A^{j,k},
$$
and the product morphism from the left side into $A$ is equal to this morphism
plus pieces which go into $\bigoplus _{j'+k'\geq j+k+1}A^{j',k'}$. This implies
the statement $(\ast )$  which we are trying to prove. Thus in the exact
sequence
of the previous paragraph, the filtration induced on
$(J^k/J^{k-1})^{\ast}$ by the weight filtration of the middle term, is equal to
the weight filtration.

The morphism $P\otimes  Sym ^{k-2} (N^{\ast}) \rightarrow Sym ^{k}(N^{\ast})$
is compatible with the filtrations induced by $W$.  We claim that $N$ and $P$
are mixed twistor structures. Thus the morphism is a morphism of mixed twistor
structures, and the cokernel is a mixed twistor structure. This cokernel
is equal to $(J^k/J^{k-1})^{\ast}$ with its filtration $W$. Thus
$Gr ^V_{k}Q_{\infty , \eta}(R_{Del}(X,G,x))$ is a mixed twistor structure so, as
explained above with the lemma of two filtrations, we are done.

We now prove that $N$ and $P$ are mixed twistor structures.
First of all, as in \cite{Goldman-Millson} \cite{Moduli},
$N$ is calculated as the
first cohomology of the complex
$$
\ker (\xi \Aa^0_X(End (E)\rightarrow End (E_x))\rightarrow
\xi \Aa ^1_X(End (E))\rightarrow \xi \Aa ^2_X(End(E)).
$$
Thus there is an exact sequence
$$
0\rightarrow H^0(End (E))\rightarrow
End (E_x) \rightarrow N \rightarrow H^1(X\times \pp ^1 / \pp ^1, End
(E))\rightarrow 0
$$
and the image of $End (E_x)$ in $N$ is the tangent space of the orbit, that
is to
say that it is $W_0$.  From \ref{coho1} the cohomology
$H^i(X\times \pp ^1 / \pp ^1, End
(E))$ is a pure twistor structure of weight $i$ (note that $End(E)$ is pure of
weight zero).  From $i=0$ and the fact that $End (E_x)$ is pure of weight
zero we find that $W_0$ is pure of weight zero. From $i=1$ we get that
$H^1(End(E))$, which  is the
piece $W_1/W_0$, is pure of weight $1$.  This proves that $N$ is a mixed twistor
structure.

Next, we have a cup-product $N \otimes N \rightarrow H^2(End (E))$ which is
actually symmetric (the antisymmetry on odd-degree classes cancels the
antisymmetry of the Lie bracket).
This gives a morphism $Sym ^2(N)\rightarrow H^2(End(E))$.  The part
$W_1Sym ^2(N)$ goes to zero since, by the group action, there are no
obstructions in the direction of the group orbit. It follows that the cokernel
and the image are
pure of weight $2$. The module of relations $P$ is the image of the dual
morphism
$$
H^2(End(E))^{\ast} \rightarrow Sym ^2(N^{\ast}),
$$
so $P$ is a pure twistor structure of weight $-2$.

This completes the proof that the jet spaces of the representation space
have natural mixed twistor structures.
\eop

{\em Application to the relative completion:}  Taking some sort of Tannakian
dual of the completions of the representation spaces $R(X,G,x)$ where $G$ runs
through all groups, at points corresponding to all representations associated to
a given semisimple representation $\rho$, should give back the {\em relative
Malcev completion} of $\pi _1(X,x)$ at the representation $\rho$. One should be
able to use this to obtain a mixed twistor structure on the relative Malcev
completion. Alternatively, a direct generalization of Hain's technique
\cite{Hain} should also work. One would like to check, in fact, that the two
constructions agree. As you can tell from the tense of this paragraph, I haven't
looked into this at all!

{\em Application to the moduli space:}
The group $G$ doesn't act on the formal local ring of $R(X,G,x)$ at a
point $\rho$ (semisimple representation), however we do obtain an infinitesimal
action of the Lie algebra ${\bf g}$. Assume that $G$ is connected. Then the
formal local ring of $M(X,G)$ at  $\rho$ is
the subring of ${\bf g}$-invariants in the  formal local ring of $R(X,G,x)$ at
$\rho$. This latter has a mixed twistor structure (it is the dual of
$J^{\infty}$). The action of ${\bf g}$ preserves the mixed twistor
structure so the
subring of invariants is a mixed twistor structure. Again
taking the dual, $J^{\infty}_{\rho}(M(X,G))$ has a mixed twistor structure. Note
that if $\rho$ is not a smooth point, the weight filtration is not at all the
same as the easy filtration, and $J^1$ might very well start of with pieces of
relatively high weight (essentially because the lower-degree terms don't survive
as invariants).

\end{document}